\documentclass{article}
\usepackage{amsmath}
\usepackage{amssymb}
\usepackage{amsthm}
\usepackage{mathrsfs}
\usepackage{graphicx}
\usepackage{subfigure}
\usepackage{float}
\usepackage{setspace}
\usepackage{extarrows}
\usepackage{hyperref}
\usepackage{bm}
\usepackage[a4paper,centering]{geometry}
\hypersetup{colorlinks=true,linkcolor=blue,anchorcolor=blue,citecolor=blue}
\numberwithin{equation}{section}
\newtheorem{definition}{\textbf{Definition}}[section]
\newtheorem{lemma}[definition]{\textbf{Lemma}}
\newtheorem{corollary}[definition]{\textbf{Corollary}}
\newtheorem{theorem}[definition]{\textbf{Theorem}}
\newtheorem{proposition}[definition]{\textbf{Proposition}}

\allowdisplaybreaks[0]
\title{The intermediate level-sets of the four-dimensional membrane model}
\author{Xinyi Li\thanks{Beijing International Center for Mathematical Research, Peking University. Research supported by the National Key R\&D Program of China (No.\ 2021YFA1002700 and No.\ 2020YFA0712900) and NSFC (No.\ 12071012). }  \and Runsheng Liu\thanks{School of Mathematical Sciences, Peking University. Research partially supported by the elite undergraduate training program of School of Mathematical Sciences in Peking University.}}
\begin{document}
\maketitle
\begin{abstract}
In this paper, we consider the discrete membrane model in four dimensions. We confirm the existence of the scaling limit of the intermediate (i.e., a multiple of the expected maximum) level-sets of the model, and show that it is equal in law to a tilted version of the sub-critical Gaussian multiplicative chaos (GMC) measure of the continuum membrane model.
\end{abstract}
\section{Introduction and main results}
The membrane model, also known as the discrete bi-Laplacian field, is a Gaussian lattice model whose Hamiltonian is given through the bi-Laplacian operator. It gets its name from counterpart continuum models in mechanics that describe the behaviour of thin structures (e.g.\ membranes) that react elastically to external forces. In \cite{Sakagawa} Sakagawa investigates the entropic repulsion phenomenon for a class of finite-range Gaussian fields, initiating the study of the membrane model within the probability community. The study of the extremal theory of the membrane model is largely inspired by that of the planar Gaussian free field (GFF) -- e.g., maximum and entropic repulsion in \cite{ref19,Daviaud}, intermediate level-sets in \cite{ref1}, extremal process in \cite{ref13,ref5});  see \cite{ref7} for a more thorough review. Kurt \cite{ref3,Kurt2} studies the maximum and the entropic repulsion of the membrane model in its critical dimension (4D). Schweiger \cite{ref2} analyzes the asymptotics of the Green function corresponding to the membrane model on the cubic lattice. An object closely related to extrema is the set of high points, a.k.a.\ level sets. In \cite{ref12} Cipriani studies the high points of the membrane model in four dimensions. There is also an on-going work \cite{GLZ} by the first author and collaborators investigating the extremal process in four dimensions, aiming at establishing similar characterizations as those for the planar GFF. In supercritical dimensions, \cite{ref16} studies the percolation of the level-sets, inspired by \cite{ref15} for similar results on supercritical GFF. The membrane model also has a scaling limit, which is introduced and studied in \cite{ref14} by Cipriani, Dan and Hazra.
\par
We now describe our model and results in detail. The discrete membrane model (with Dirichlet boundary conditions) on a finite subset $V\subset\mathbb{Z}^{d}$ is a centered Gaussian field $\{h_{x}\}$, which vanishes outside $V$, whose Hamiltonian is given by $\frac{1}{2}\sum_{x\in V}|\Delta h_{x}|^{2}$, where $\Delta$ denotes the discrete Laplacian operator. We will give its precise definition in the next section.
\par
In 2009, Kurt \cite{ref3} shows that the extreme value of the 4D membrane model in square lattices $D_{N}:=(0,N)^{4}\cap\mathbb{Z}^{4}$ grows as
\[\max_{x\in D_{N}}h_{x}=2\sqrt{2\gamma}\log N(1+o(1)),\]
where $o(1)\to0$ in probability as $N\to\infty$ and $\gamma:=8/\pi^{2}$ satisfies that the asymptotic of the Green function with respect to the 4D membrane model (see Definition \ref{gf}) satisfies $G^{D_{N}}(x,x)=\gamma\log N+O(1)$ as $N\to\infty$ for $x$ ``deep'' inside $D_{N}$. \cite{ref3} also discussed the intermediate level-sets
\[\{x\in D_{N}:h_{x}\geq2\lambda\sqrt{2\gamma}\log N\},\quad\lambda\in(0,1),\]
and shows the size of the intermediate level-sets is $N^{4(1-\lambda^{2})+o(1)}$, where $o(1)\to0$ in probability as $N\to\infty$.
\par
In this paper, we are going to show that the point measure induced by the intermediate level-sets, namely
\[\sum_{x\in D_{N}}\delta_{x/N}\otimes\delta_{h_{x}^{D_{N}}-a_{N}},\]
with $a_{N}$ a scale sequence such that for some $\lambda\in(0,1)$,
\begin{equation}\label{000}\lim_{N\to\infty}\frac{a_{N}}{\log N}=2\lambda\sqrt{2\gamma},\end{equation}
has a scaling limit which is non-trivial and can be characterized in an explicit way. More precisely, setting scaling sequence
\begin{equation}\label{111}K_{N}:=\frac{N^{4}}{\sqrt{\log N}}e^{-\frac{a_{N}^{2}}{2\gamma\log N}},\end{equation}
we have:
\begin{theorem}\label{2.1}
Set $D=(0,1)^{4}$ and $D_{N}=(0,N)^{4}\cap\mathbb{Z}^{4}$. For each $\lambda\in(0,1)$ and each sample $h^{D_{N}}$ of the membrane model in $D_{N}$, define the rescaled point measure
\begin{equation}\label{112}\eta_{N}^{D}:=\frac{1}{K_{N}}\sum_{x\in D_{N}}\delta_{x/N}\otimes\delta_{h_{x}^{D_{N}}-a_{N}},\end{equation}
with $a_{N},K_{N}$ as in \eqref{000} and \eqref{111}. Then, there is a random Borel measure $Z_{\lambda}^{D}$ on $\overline{D}$ with $\mathbb{E}[Z_{\lambda}^{D}(\overline{D})]\in\mathbb{R}_{+}$ such that
\begin{equation}\label{113}\eta_{N}^{D}\xrightarrow[N\to\infty]{\mathrm{law}}Z_{\lambda}^{D}(\mathrm{d}x)\otimes e^{-\pi\lambda h}\mathrm{d}h\end{equation}
with respect to the topology of vague convergence of measures on $\overline{D}\times\mathbb{R}$.
\end{theorem}
As a direct consequence of Theorem \ref{2.1}, we can give a sharp asymptotic for the size of the intermediate level-sets:
\begin{corollary}\label{2.2}
For $a_{N},K_{N}$ as in \eqref{000} and \eqref{111},
\[\frac{1}{K_{N}}\#\{x\in D_{N}:h^{D_{N}}(x)\geq a_{N}\}\xrightarrow[N\to\infty]{\mathrm{law}}(\pi\lambda)^{-1}Z_{\lambda}^{D}(\overline{D}).\]
\end{corollary}
In fact, the measure $Z_{\lambda}^{D}$ can be characterized by the scaling limit of membrane model. In \cite{ref14}, Cipriani introduces the continuum membrane model which serves as a counterpart of the membrane model in the continuum. The following theorem shows that the measure $Z_{\lambda}^{D}$ is indeed a tilted version of the Gaussian multiplicative chaos (GMC) measure corresponding to the continuum membrane model on $D$:
\begin{theorem}\label{2.4}
For $D=(0,1)^{4}$, the random measure $Z_{\lambda}^{D}$ has the following law:
\[Z_{\lambda}^{D}(\mathrm{d}x)\overset{\mathrm{law}}{=}\frac{\sqrt{\pi}}{4}e^{\frac{4\lambda^{2}}{\gamma}s_{D}(x)}\mu_{\infty}^{D,\pi\lambda}(\mathrm{d}x)\]
where the $\mu_{\infty}^{D,\pi\lambda}$ is the subcritical GMC measure on $D$ defined in \eqref{gmc} and the function $s_{D}(x)$ is defined in \eqref{sdx}.
\end{theorem}
The proof of Theorem \ref{2.1} is inspired by \cite{ref1}, which discussed scaling limits of intermediate level-sets of DGFF. We first prove that $\eta_{N}^{D}$ is tight with respect to the vague topology via a truncated second moment argument, then show that the subsequential limits of $\eta_{N}^{D}$ can be factorized as the right-hand side of \eqref{113}, and finally we prove the uniqueness of the subsequential limits by showing that they must satisfy some axiomatic characterization (see Proposition \ref{3.8} and \ref{3.11}). Theorem \ref{2.4} is an application of Shamov's approach to Gaussian Multiplicative Chaos (see \cite{ref11} for more details).
\par
In principle, the results in Theorems \ref{2.1}-\ref{2.4} may be also extended to general nice domains in $\mathbb{R}^{4}$. However, the analysis of the asymptotics of the Green function corresponding to the discrete bi-Laplacian operator is much more difficult than the case of GFF, for which the random walk representation plays a crucial role. Although a representation by intersections of two independent random walks exists for the discrete bi-Laplacian operator, it does not satisfy the boundary condition of our model. Deriving these asymptotics falls outside the scope of this work, therefore we do not pursue this generalization any further.
\par
Another possible direction for generalization is to extend the results above to the Gaussian polylaplacian models (in their critical dimension), which are defined by replacing the discrete bi-Laplacian in \eqref{dmm} by $\nabla\Delta^{\frac{n-1}{2}}$ or $\Delta^{\frac{n}{2}}$ depending on whether the parameter $n$ (which also indicates the respective critical dimension) is odd or even. The key points in the proof are the Gibbs-Markov property and the asymptotics of the Green function. The Gibbs-Markov property holds for all such models. The $n$-dimensional Gaussian poly-Laplacian model with even parameter $n$ should be log-correlated by applying the methods in \cite{ref8,ref9,ref4}. Then the similar scaling limit results of the intermediate level-sets should follow when one derives sufficiently sharp asymptotics of the Green function. When $n$ is odd, since the operator is no longer local, the scheme above does not work any more. It is an interesting question whether one can still obtain similar results in this case.
\par
We now briefly describe how this paper is organized. In Section 2, we will give some brief introduction on the discrete and continuum membrane model, and clarify the notation in this paper. In Section 3, we will estimate the first moment and the second moment of the size of truncated level-sets. In Section 4, we will prove the existence and factorization property of the subsequential limits. In Section 5, we will characterize the limit measure, and prove the uniqueness of the subsequential limits, thus complete the proof of Theorem \ref{2.1}. In Section 6, we will prove Theorem \ref{2.4}.

\section{Notation and Preliminaries}
In this section, we will introduce our notation, and then give some basic results on the membrane model.
\subsection{Notation}
In this paper, $\nabla,\Delta$ denotes the discrete gradient and Laplacian operator, which is defined by:
\[\nabla f(x):=(f(x+e_{1})-f(x),\cdots,f(x+e_{d})-f(x)),\]
for $e_{i}$ the $i$-th unit coordinate vector in $\mathbb{R}^{d}$, and
\[\Delta f(x)=\frac{1}{2d}\sum_{y:y\sim x}(f(y)-f(x)),\]
where the sum is taken over the neighbors of $x$. Let $\overline{\nabla}$, (resp.\ $\widetilde{\Delta}$) denote the continuum gradient (resp.\ Laplacian) operator. Let $G,\widetilde{G}$ denote the discrete and continuum Green function, which will be defined in the next subsection. For a domain $D$, let $H_{0}^{2}(D)$ be the Sobolev space of index 2 on $D$, that is, the completion of the test function space in $D$ with respect to the inner product $\langle\cdot,\cdot\rangle_{\widetilde{\Delta}}$, where
\[\langle f,g\rangle_{\widetilde{\Delta}}:=\int_{D}\widetilde{\Delta} f\widetilde{\Delta} g.\]
\par
The constants in this paper, such as $c,c',\widetilde{c},\cdots$, are universal, but may change from place to place. And $c(r,s,\cdots)$ denotes a constant depending only on $r,s,\cdots$, whose precise value may change from place to place. Moreover, these constants are assumed to be positive without further specification.
\par
Let $\mathfrak{D}$ be the collection of all open dyadic cubes and their finite unions. The dyadic cubes takes the form
\[(i2^{-n},(i+1)2^{-n})\times(j2^{-n},(j+1)2^{-n})\times(k2^{-n},(k+1)2^{-n})\times(l2^{-n},(l+1)2^{-n}),\quad i,j,k,l,n\in\mathbb{Z}.\]
For $D\in\mathfrak{D}$, let $D_{N}$ denote its lattice approximation, which is a subset of $\mathbb{Z}^{4}$, such that $D_{N}/N\subset D$ and for any $x\in D^{c}$, the $l^{\infty}$ distance between $Nx$ and $D_{N}$ is larger than 1. One can assume that every domain $D$ in this paper belongs to $\mathfrak{D}$.
\par
Let $\|\cdot\|$ be the $l^{\infty}$ norm of a vector.
\subsection{Preliminaries}
We start with the definition of the (discrete) membrane model:
\begin{definition}
The discrete membrane model (with Dirichlet boundary condition) on a finite subset $V$ of $\mathbb{Z}^{d}$ is a centered Gaussian field denoting by $h_{x}^{V}$, vanishes outside $V$, and has the probability density function
\begin{equation}\label{dmm}
\mathbb{P}(\mathrm{d}h^{V})\propto\exp\left(-\frac{1}{2}\sum_{x\in\mathbb{Z}^{d}}|\Delta h_{x}|^{2}\right)\prod_{x\in V}\mathrm{d}h_{x}\prod_{x\notin V}\delta_{0}(\mathrm{d}h_{x}),
\end{equation}
where $\delta_{0}$ is the Dirac point mass at 0.
\end{definition}
Next, we will define the Green function of the membrane model.
\begin{definition}\label{gf}
The Green function $G^{V}(x,y)$ is the correlation function of the membrane model, that is, $G^{V}(x,y):=\mathbb{E}[h_{x}^{V}h_{y}^{V}]$.
\end{definition}
The Green function solves the following difference equation:
\[\begin{cases}\Delta^{2}G^{V}(x,y)=\delta_{x}(y),\ \ & y\in V,\\ G^{V}(x,y)=0,\ \ & y\in\partial_{2}V,\end{cases}\]
where $\partial_{2}V$ denotes the points in $V^{c}$ whose graph distance to $V$ is less or equal than 2.
\par
The estimation of the Green function is important throughout the proof. We now paraphrase the asymptotics from \cite[Observation 1.5]{ref2}. Before we state it, we need to define the continuum Green function $\widetilde{G}$:
\begin{definition}
The continuum Green function of bi-Laplacian operator on a domain $D$ with Dirichlet boundary condition on $\partial D$ is a function $\widetilde{G}^{D}(x,y)=\widetilde{G}_{y}^{D}(x)$ that solves the following equations:
\begin{equation}\label{cgf}\begin{cases}\widetilde{\Delta}^{2}\widetilde{G}_{y}(x)=\delta_{y}(x),\ \ & x\in D,\\ \widetilde{G}_{y}(x)=0,\ \ & x\in\partial D,\\ \frac{\partial\widetilde{G}_{y}}{\partial\bm{n}}(x)=0,\ \ & x\in\partial D.\end{cases}\end{equation}
\end{definition}
Recall $D=(0,1)^{4}$, $D_{N}=(0,N)^{4}\cap\mathbb{Z}^{4}$ and $\gamma=8/\pi^{2}$. Letting $d(x)$ be the $l^{\infty}$ distance from $x$ to $\partial D$, denoting $\Gamma(x,y):=\gamma\log|x-y|$ as the fundamental solution, setting
\begin{equation}\label{adxy}a^{D}(x,y):=\widetilde{G}^{D}(x,y)-\Gamma(x,y)\quad \forall x,y\in D\ \mathrm{and}\ x\neq y,\end{equation}
and thanks to \cite[Lemma 3.5]{ref2}, $a^{D}(x,y)$ can be continuously extended to $D\times D$, and we define
\begin{equation}\label{sdx}s_{D}(x):=a^{D}(x,x).\end{equation}
We then have:
\begin{theorem}[{\cite[Observation 1.5]{ref2}}]\label{gfa}
There exists $\theta>0$, such that:
\begin{equation}\label{251}G^{D_{N}}([xN],[xN])=\gamma\log N+s_{D}(x)+o(1)\end{equation}
with $o(1)\to0$ as $N\to\infty$ for any $x\in D$ with $d(x)\geq\frac{(\log N)^{\theta}}{N}$. And
\begin{equation}\label{252}G^{D_{N}}([xN],[yN])=\widetilde{G}^{D}(x,y)+o(1)\end{equation}
with $o(1)\to0$ as $N\to\infty$ for all $x,y\in\{(x,y):x,y\in D,x\neq y\}$ with $\min\{d(x),d(y)\}\geq\frac{(\log N)^{\theta}}{N}$. Moreover
\begin{equation}\label{253}\left|G^{D_{N}}(x,y)-\gamma\log\left(2+N\frac{\max(d(x/N),d(y/N))}{1+|x-y|}\right)\right|=O(1).\end{equation}
\end{theorem}
Since the above asymptotics are scale-invariant, and the Green function of the union of disjoint sets is as simple as the Green function on each set. Therefore the above asymptotics \eqref{251}-\eqref{253} are valid on $\mathfrak{D}$.
\par
Similar to the Gaussian free field, the membrane model also enjoys the Gibbs-Markov property. We formalize it in the following theorem and its proof will be given in the Appendix.
\begin{theorem}\label{gmp}
For $U\subset V\subset\mathbb{Z}^{d}$, let $h^{V}$ be the membrane model in $V$. Define
\[\varphi^{V,U}(x):=\mathbb{E}[h_{x}^{V}|\sigma(h_{z}^{V}:z\in V\backslash U)],\quad x\in U.\]
Then $\{h^{V}-\varphi^{V,U}\}_{x\in U}$ is equal in law to a discrete membrane model in $U$, and is independent of $\varphi^{V,U}$.
\end{theorem}
By Theorem \ref{gmp}, we have $G^{V}(x,x)-G^{U}(x,x)=\mathrm{var}[\varphi^{V,U}(x)]$, which gives
\begin{equation}\label{263}G^{V}(x,x)\geq G^{U}(x,x),\quad x\in U\subset V.\end{equation}
The continuum membrane model, introduced in \cite{ref14} can be seen as the scaling limit of the discrete membrane model. We now give two equivalent definitions of the continuum membrane model.
\begin{definition}\label{cmm}
A continuum membrane model on a bounded, open $D\subset\mathbb{R}^{4}$ is a functional $h$ that maps each bounded, measurable $f:D\to\mathbb{R}$ to a random variable such that
\par
(1) $h$ is a.s. linear, i.e.,
\[h(af+bg)=ah(f)+bh(g)\quad\mathrm{a.s.}\]
for any bounded measurable $f,g$ and any $a,b\in\mathbb{R}$, and
\par
(2) for any bounded measurable $f$, $h(f)\overset{\mathrm{law}}{=}\mathcal{N}(0,\sigma^{2}(f))$, where
\[\sigma^{2}(f)=\int_{D\times D}\widetilde{G}^{D}(x,y)f(x)f(y)\mathrm{d}x\mathrm{d}y.\]
Equivalently, $h$ is a centered Gaussian process on the space of bounded, measurable function $f:D\to\mathbb{R}$ such that the covariance structure is given by:
\[\mathrm{cov}(h(f),h(g))=\int_{D\times D}\widetilde{G}^{D}(x,y)f(x)g(y)\mathrm{d}x\mathrm{d}y.\]
\end{definition}
Alternatively, the continuum membrane model can be defined as the canonical Gaussian process on Hilbert space $H_{0}^{2}(D)$ with the inner product $\langle f,g\rangle_{\widetilde{\Delta}}=\int_{D}\widetilde{\Delta}f(x)\widetilde{\Delta}g(x)\mathrm{d}x$. If we let $\{f_{k}\}$ be an orthonormal basis of the Hilbert space, then with $Z_{k}$ i.i.d. standard normal random variables, we will denote the partial sums by
\[\phi_{n}(x):=\sum_{k=1}^{n}Z_{k}f_{k}(x).\]
For any smooth function $f\in H_{0}^{2}(D)$, let $h_{n}(f):=\int_{D}f(x)\phi_{n}(x)\mathrm{d}x$, then the $L^{2}$-limit of $h_{n}$ in $H_{0}^{2}(D)$ satisfies Definition \ref{cmm}.
\par
Finally, we will define the Gaussian multiplicative chaos associated with the continuum membrane model. For each $\beta\in\mathbb{R}_{+}$, define the random measure
\[\mu_{n}^{D,\beta}(\mathrm{d}x):=1_{D}(x)e^{\beta\phi_{n}(x)-\frac{1}{2}\beta^{2}\mathbb{E}[\phi_{n}(x)^{2}]}\mathrm{d}x.\]
Thanks to Kahane \cite{ref10} (see also Shamov's approach \cite{ref11}), there exists an a.s.\ finite random Borel measure $\mu_{\infty}^{D,\beta}$, supported on $D$, and for each measurable $A\subset D$,
\begin{equation}\label{gmc}\lim_{n\to\infty}\mu_{n}^{D,\beta}(A)=\mu_{\infty}^{D,\beta}(A),\quad\mathrm{a.s.}\end{equation}
Then the measure $\mu_{\infty}^{D,\beta}$ is defined as the GMC associated with the continuum membrane model, and such measure is unique regardless of the choice of the orthonormal basis.
\par
Finally, we introduce a ``reversed Jensen'' inequality from \cite[Lemma 3.12]{ref1}.
\begin{lemma}\label{lem:rev}
Let $X_{1},\ldots,X_{n}$ be non-negative random variables. For any $\varepsilon>0$, we have
\[\mathbb{E}[\exp\{-(X_{1}+\cdots+X_{n})\}]\leq\exp\left\{-e^{-\varepsilon}\sum_{i=1}^{n}\mathbb{E}[X_{i}1_{\{X_{i}\leq\varepsilon\}}]\right\}\]
\end{lemma}
\section{Moments of the size of level-sets}
We will start our proof by calculating the first and the second moments of the size of level-sets. Lemma \ref{3.1} and \ref{3.2} calculate the first moments of the size of level-sets, and when we turn to the second moment estimations, a truncation is needed for all $\lambda\in(0,1)$. Lemma \ref{4.1} discusses the difference between the truncated and untruncated level-sets, and Proposition \ref{4.2} calculates the truncated second moments.
\par
For each $b\in\mathbb{R}$, define
\[\Gamma_{N}^{D}(b):=\{x\in D_{N}:h^{D_{N}}(x)\geq a_{N}+b\}.\]
First of all, we will give a general size bound of $\Gamma_{N}^{D}(b)$:
\par
\begin{lemma}\label{3.1}
For any $D\in\mathfrak{D}$ there is a $c=c(D)$ such that for any $N\geq1$, any $b\in[-\log N,\log N]$, any sequences $a_{N}$ that satisfies \eqref{000} and any $A\subset D_{N}$, we have
\begin{equation}\label{311}\mathbb{E}|\Gamma_{N}^{D}(b)\cap A|\leq cK_{N}\frac{|A|}{N^{4}}e^{-\frac{a_{N}b}{\gamma\log N}}.\end{equation}
\end{lemma}
\begin{proof}
We will prove \eqref{311} by summing over $x\in D_{N}$. In fact, it suffices to prove:
\[\mathbb{P}(h^{D_{N}}(x)\geq a_{N}+b)\leq\frac{c(D)}{\sqrt{\log N}}e^{-\frac{a_{N}^{2}}{2\gamma\log N}}e^{-\frac{a_{N}b}{\gamma\log N}}\]
uniformly in $x\in D_{N}$ and $|b|\leq\log N$. Without loss of generality, we can assume that $D$ contains only one dyadic square.
\par
By Theorem \ref{gmp}, for lattice domains $U\subset V\subset\mathbb{Z}^{d}$, we can write $h^{V}\overset{\mathrm{law}}{=}\varphi^{V,U}+h^{U}$, and the two parts in the right-hand side are independent. Then we have
\begin{align*}
\mathbb{P}(h^{U}(x)\geq a)&\leq\mathbb{P}(h^{U}(x)\geq a,\varphi^{V,U}(x)\geq0)+\mathbb{P}(h^{U}(x)\geq a,\varphi^{V,U}(x)\leq0)\\
&=2\mathbb{P}(h^{U}(x)\geq a,\varphi^{V,U}(x)\geq0)\leq2\mathbb{P}(h^{V}(x)\geq a).
\end{align*}
We then enlarge $D_{N}$ to a larger lattice $\widetilde{D}_{N}$, whose diameter is 2 times larger than that of $D_{N}$, and ensures that $\widetilde{D}_{N}$ and $D_{N}$ have the same center. Then the variance of $h^{\widetilde{D}_{N}}(x)(x\in D_{N})$ is uniformly within a constant of $\gamma\log N$ thanks to \eqref{253}. So for some $c'$,
\begin{align*}
\mathbb{P}(h^{\widetilde{D}_{N}}(x)\geq a_{N}+b)&\leq\frac{1}{\sqrt{2\pi}}\frac{1}{\sqrt{\gamma\log N-c'}}\int_{b}^{\infty}e^{-\frac{(a_{N}+s)^{2}}{2(\gamma\log N+c')}}\mathrm{d}s\\
&\leq\frac{1}{\sqrt{2\pi}}\frac{1}{\sqrt{\gamma\log N-c'}}\int_{b}^{\infty}e^{-\frac{a_{N}^{2}+2a_{N}s}{2(\gamma\log N+c')}}\mathrm{d}s\\
&\leq\frac{C}{\sqrt{\log N}}e^{-\frac{a_{N}^{2}}{2\gamma\log N}}e^{-\frac{a_{N}b}{\gamma\log N}}.
\end{align*}
This completes the proof of the lemma.
\end{proof}
\par
Thanks to Lemma \ref{3.1}, more accurate asymtotics of the above expectation can be calculated for open sets $A$ whose closure lies inside $D$.
\begin{lemma}\label{3.2}
For any $b\in\mathbb{R}$ and any open sets $A$ whose closure lie inside $D$,
\[\mathbb{E}|\{x\in\Gamma_{N}^{D}(b):x/N\in A\}|=\frac{e^{-\pi\lambda b}}{4\lambda\sqrt{\pi}}\left(\int_{A}e^{\frac{4\lambda^{2}}{\gamma}s_{D}(x)}\mathrm{d}x+o(1)\right)K_{N},\]
with $o(1)\to0$ as $N\to\infty$ uniformly on compact sets of $b$.
\end{lemma}
\begin{proof}
Since the closure of $A$ lies in $D$, all the points in $A$ lies ``deep'' inside $D$. Thanks to \eqref{251} we have
\[G^{D_{N}}([xN],[xN])=\gamma\log N+s_{D}(x)+o(1),\]
with $o(1)\to0$ as $N\to\infty$ uniformly in $x\in\overline{A}$. Then using the following Gaussian tail bound approximation from \cite[Theorem 1.2.6]{ref17}:
\begin{equation}\label{321}\left(\frac{1}{x}-\frac{1}{x^{3}}\right)\frac{e^{-\frac{x^{2}}{2}}}{\sqrt{2\pi}}\leq\mathbb{P}(X\geq\sigma x)\leq\frac{1}{x}\frac{e^{-\frac{x^{2}}{2}}}{\sqrt{2\pi}},\end{equation}
with $X\sim N(0,\sigma^{2})$. We have:
\begin{align*}
\mathbb{P}(h^{D_{N}}(x)\geq a_{N}+b)&=(1+o(1))\frac{1}{\sqrt{2\pi}}\frac{\sqrt{\gamma\log N+s_{D}(x)+o(1)}}{a_{N}+b}e^{-\frac{(a_{N}+b)^{2}}{2(\gamma\log N+s_{D}(x)+o(1))^{2}}}\\
&=(1+o(1))\frac{1}{4\lambda\sqrt{\pi}\sqrt{\log N}}e^{-\frac{a_{N}^{2}}{2\gamma\log N}}e^{-\pi\lambda b}e^{\frac{4\lambda^{2}s_{D}(x)}{\gamma}}.
\end{align*}
Since $s_{D}(x)$ is continuous in $D$, we can replace the Riemann sum by integrating on $A$, thus completes the proof.
\end{proof}
\par
Then, we need to deal with a truncated version of the measures $\eta_{N}^{D}$ for the proof for all $\lambda\in(0,1)$. Let $\{D_{N}\}$ be the lattice approximation of a domain $D\in\mathfrak{D}$. Denote $\Lambda_{r}(x):=\{z\in\mathbb{Z}^{4}:\|z-x\|\leq r\}$ and, for each $N\geq1$ and each $x\in D_{N}$, let
\[n(x):=\max\{n\geq0:\Lambda_{e^{n+1}}(x)\subset D_{N}\}.\]
Thanks to Theorem \ref{gfa}, we have $\log N-c(\varepsilon)\leq n(x)\leq\log N+c'(D)$ for all $x\in D_{N}$ such that $d(x,D_{N}^{c})>\varepsilon N$. We now define the sequence of domains
\[\Delta^{k}(x):=\begin{cases}\varnothing,\ \ & k=0,\\ \Lambda_{e^{k}}(x),\ \ & 1\leq k\leq n(x)-1,\\ D_{N},\ \ & k=n(x).\end{cases}\]
For $U\subset V$, recall $\varphi^{V,U}$ as in Theorem \ref{gmp}. We then set
\[S_{k}(x):=\varphi^{D_{N},\Delta^{k}(x)}(x),\quad k=0,1,\cdots,n(x).\]
By definition, we have $S_{0}(x)=h_{x}^{D_{N}}$ and $S_{n(x)}(x)=0$.
\par
Next, for a sequence $a_{N}$ satisfying \eqref{000} and $M>0$, we define the truncation event
\[T_{N,M}(x):=\bigcap_{k=k_{N}}^{n(x)}\left\{\left|S_{k}(x)-a_{N}\frac{n(x)-k}{n(x)}\right|\leq M(n(x)-k)^{3/4}\right\}\]
where
\[k_{N}:=\left[\frac{1}{8}\log K_{N}\right]=\frac{1}{2}(1-\lambda^{2}+o(1))\log N.\]
Then we can define the truncated point measure
\begin{equation}\label{tpm}\widehat{\eta}_{N}^{D,M}:=\frac{1}{K_{N}}\sum_{x\in D_{N}}1_{T_{N,M}(x)}\delta_{x/N}\otimes\delta_{h^{D_{N}}(x)-a_{N}}.\end{equation}
By the definition of $\widehat{\eta}_{N}^{D,M}$, we immediately get $\langle\widehat{\eta}_{N}^{D,M},f\rangle\leq\langle\eta_{N}^{D},f\rangle$ for any measurable $f\geq0$.
\par
First of all, we need to deal with the difference between $\widehat{\eta}_{N}^{D,M}$ and $\eta_{N}^{D}$ when $M$ is sufficiently large. For this, we define the truncated level-sets
\[\widehat{\Gamma}_{N}^{D,M}(b):=\{x\in D_{N}:h^{D_{N}}(x)\geq a_{N}+b, T_{N,M}(x)\}.\]
\begin{proposition}\label{4.1}
For each $\lambda\in(0,1)$ and each $b_{0}>0$ there exist $c,\widetilde{c}$ such that for any $D\in\mathfrak{D}$, any $b\in[-b_{0},b_{0}]$, any $M\geq1$ and large $N$,
\[\mathbb{E}|\Gamma_{N}^{D}(b)\backslash\widehat{\Gamma}_{N}^{D,M}(b)|\leq ce^{-\widetilde{c}M^{2}}(\mathrm{diam}D)^{4+4\lambda^{2}}K_{N}.\]
\end{proposition}
We postpone its proof till the end of this section.
\par
A direct consequence of Proposition \ref{4.1} is that for any bounded, measurable $f:D\times\mathbb{R}\to\mathbb{R}$,
\[\lim_{M\to\infty}\limsup_{N\to\infty}|\langle\widehat{\eta}_{N}^{D,M},f\rangle-\langle\eta_{N}^{D},f\rangle|=0,\quad\mathrm{a.s.}\]
By Lemma \ref{3.1} and the domination of $\widehat{\eta}_{N}^{D,M}$ by $\eta_{N}^{D}$, the family of truncated measures $\{\widehat{\eta}_{N}^{D,M}:N\geq1\}$ is tight with respect to the vague topology. Then a subsequential weak limit $\widehat{\eta}^{D,M}$ exists and we can study its properties, which will be discussed in the next section.
\par
We then need to calculate the truncated second moments. For $b,b'\in\mathbb{R}$ with $b<b'$, let
\[\widehat{\Gamma}_{N}^{D,M}(b,b'):=\widehat{\Gamma}_{N}^{D,M}(b)\backslash\widehat{\Gamma}_{N}^{D,M}(b').\]
Then we have:
\begin{proposition}\label{4.2}
Let $\lambda\in(0,1)$. For any $\varepsilon,M>0$ and any $b,b'\in\mathbb{R}$ with $b<b'$, there is $c=c(\varepsilon,M,b,b')$ such that for any $D\in\mathfrak{D}$ and any large $N$,
\[\mathbb{E}[|\widehat{\Gamma}_{N}^{D,M}(b,b')\cap D_{N}^{\varepsilon}|^{2}]\leq c(\mathrm{diam}D)^{8+8\lambda^{2}}K_{N}^{2},\]
where $D_{N}^{\varepsilon}$ consists of the points in $D_{N}$ that have distance at least $\varepsilon N$ to the complement of $D_{N}$.
\end{proposition}
In order to prove Proposition \ref{4.1} and \ref{4.2}, we need the following lemmas:
\begin{lemma}\label{4.6}
For any $\varepsilon,r>0$, there is $c=c(\varepsilon,r)$ such that for any $D\in\mathfrak{D}$ with $\mathrm{diam}D\leq r$ and any large $N$, we have:
\par
(1) For any $x\in D_{N}$ and any $k_{N}\leq k\leq m\leq n(x)$,
\begin{equation}\label{351}\mathrm{var}[S_{k}(x)-S_{m}(x)]=(m-k)\gamma+o(1)\end{equation}
where $o(1)\to0$ when $N\to\infty$.
\par
(2) For any $x\in D_{N}^{\varepsilon}$ and any $k$ such that $k_{N}\leq k\leq n(x)$,
\begin{equation}\label{352}\mathrm{var}[S_{k}(x)]-(n(x)-k)\gamma\in[0,c].\end{equation}
(3) For any $l\geq1$ there is $c'=c'(\varepsilon,l)$ such that for any $x\in D_{N}^{\varepsilon}$, any $k$ such that $k_{N}\leq k\leq n(x)$, any $m$ such that $k-l\leq m\leq k$ and any $y\in D_{N}$ satisfies $\Delta^{m+1}(y)\subseteq\Delta^{k}(x)\backslash\{x\}$, we have
\begin{equation}\label{353}|\mathbb{E}[S_{k}(x)S_{m}(y)]-(n(x)-k)\gamma|\leq c',\quad\mathrm{var}[S_{m}(y)-S_{k}(x)]\in[\gamma/2,c'].\end{equation}
\end{lemma}
\begin{proof}
By \eqref{251} and Theorem \ref{gmp},
\[\mathrm{var}[S_{k}(x)-S_{m}(x)]=\mathrm{var}\left[\varphi^{\Lambda_{e^{m}}(0),\Lambda_{e^{k}}(0)}(0)\right]=G^{\Lambda_{e^{m}}(0)}(0,0)-G^{\Lambda_{e^{k}}(0)}(0,0)=(m-k)\gamma+o(1),\]
which gives \eqref{351}. For \eqref{352}, if $k=n(x)$ the claim is automatically true, thus we only need to consider the case when $k<n(x)$. Then, since $n(x)\geq\log N-c$ for $c=c(\varepsilon)>0$, we may find $\widetilde{c}=\widetilde{c}(\varepsilon,r)>0$ such that $D_{N}\subseteq\Lambda_{e^{n(x)+\widetilde{c}}}(x)$. Thanks to \eqref{263}, we have
\[\mathrm{var}[S_{k}(x)]\leq G^{\Lambda_{e^{n(x)+\widetilde{c}}}(0)}(0,0)-G^{\Lambda_{e^{k}}(0)}(0,0)\leq(n(x)+\widetilde{c}-k)\gamma+o(1).\]
By definition $D_{N}\supseteq\Lambda_{e^{n(x)+1}}(x)$, we thus have
\[\mathrm{var}[S_{k}(x)]\geq G^{\Lambda_{e^{n(x)+1}}(0)}(0,0)-G^{\Lambda_{e^{k}}(0)}(0,0)\geq(n(x)+1-k)\gamma+o(1).\]
This completes the proof of \eqref{352}.
\par
For \eqref{353}, note that $\mathbb{E}[S_{k}(x)S_{m}(y)]=\mathbb{E}[\varphi^{D_{N},\Delta^{k}(x)}(x)\varphi^{D_{N},\Delta^{k}(x)}(y)]$ since Theorem \ref{gmp} permits us to write $\varphi^{U,W}\overset{\mathrm{law}}{=}\varphi^{U,V}+\varphi^{V,W}$ for $W\subset V\subset U$. Then, by Theorems \ref{gfa} and \ref{gmp}, this expectation equals
\[\mathbb{E}[h^{D_{N}}(x)h^{D_{N}}(y)]-\mathbb{E}[h^{\Delta^{k}(x)}(x)h^{\Delta^{k}(x)}(y)]\in[(n(x)-k)\gamma-c',(n(x)-k)\gamma+c']\]
where $c'>0$ depends on the distance of $y$ to the boundary of $\Delta^{k}(x)$, and thus the choice of $l$.
\par
Finally, we have
\[\mathrm{var}[S_{m}(y)-S_{k}(x)]=\mathrm{var}[S_{m}(y)]+\mathrm{var}[S_{k}(x)]-2\mathbb{E}[S_{m}(y)S_{k}(x)],\]
together with \eqref{352} and the first equation of \eqref{353} we get the upper bound of the second equation of \eqref{353}. For the lower bound the second equation of \eqref{353}, we have
\begin{align*}
\mathrm{var}[S_{m}(y)-S_{k}(x)]&=\mathrm{var}[\varphi^{D_{N},\Delta^{k}(x)}(y)-\varphi^{D_{N},\Delta^{k}(x)}(x)+\varphi^{\Delta^{k}(x),\Delta^{m}(y)}(y)]\\
&\geq\mathrm{var}[\varphi^{\Delta^{k}(x),\Delta^{m}(y)}(y)]\geq\mathrm{var}[\varphi^{\Delta^{m+1}(y),\Delta^{m}(y)}(y)]\geq\gamma+o(1).
\end{align*}
This proves \eqref{353}.
\end{proof}
\begin{lemma}\label{4.7}
If $x,y\in D_{N}$ satisfy $\Delta^{k}(x)\subseteq\Delta^{m}(y)$, then for any $k'<k$ and $m'>m$, the increments $S_{k'}(x)-S_{k}(x)$ and $S_{m}(y)-S_{m'}(y)$ are independent. Particularly, for any $x\in D_{N}$, the process $\{S_{k}(x)\}_{k=0}^{n(x)}$ has independent increments.
\end{lemma}
We omit its proof since it is a trivial consequence of Theorem \ref{gmp}.
\par
Now we are ready for:
\par
\begin{proof}[Proof of Proposition \ref{4.1}.] For any $b<b'$ and $\varepsilon>0$, we have:
\begin{equation}\label{331}
\begin{split}
\mathbb{E}|\Gamma_{N}^{D}(b)\backslash\widehat{\Gamma}_{N}^{D,M}(b)|\leq\ &\mathbb{E}|\Gamma_{N}^{D}(b)\backslash D_{N}^{\varepsilon}|+\mathbb{E}|\Gamma_{N}^{D}(b')|\\
&+
\sum_{x\in D_{N}^{\varepsilon}}\sum_{k=k_{N}}^{n(x)}\mathbb{P}(h^{D_{N}}(x)-a_{N}\in[b,b'),T_{N,M}^{k}(x)),
\end{split}
\end{equation}
where $T_{N,M}^{k}(x):=\{|S_{k}(x)-a_{N}\frac{n(x)-k}{n(x)}|>M(n(x)-k)^{3/4}\}$. The first two terms on the right-hand side of \eqref{331} can be bounded by Lemma \ref{3.1} by taking $b'$ sufficiently large and $\varepsilon$ small enough. Now we turn to the estimation of the third term on the right-hand side of \eqref{331}. For given $x\in D_{N}^{\varepsilon}$ and $s\in[b,b']$, we have
\[(S_{k}(x)|h^{D_{N}}(x)=a_{N}+s)\overset{\mathrm{law}}{=}\mathcal{N}\left(\frac{\mathrm{var}[S_{k}(x)]}{\mathrm{var}[S_{0}(x)]}(a_{N}+s),\frac{\mathrm{var}[S_{k}(x)]\mathrm{var}[S_{0}(x)-S_{k}(x)]}{\mathrm{var}[S_{0}(x)]}\right).\]
By Lemma \ref{4.6}, we have:
\[\left|\frac{\mathrm{var}[S_{k}(x)]}{\mathrm{var}[S_{0}(x)]}-\frac{n(x)-k}{n(x)}\right|\leq\frac{c}{n(x)}\]
and
\[\frac{\mathrm{var}[S_{k}(x)]\mathrm{var}[S_{0}(x)-S_{k}(x)]}{\mathrm{var}[S_{0}(x)]}\leq c'(n(x)-k).\]
By \eqref{321}, we then have
\[\mathbb{P}(T_{N,M}^{k}(x)|h^{D_{N}}(x)=a_{N}+s)\leq ce^{-\widetilde{c}M^{2}(n(x)-k)^{1/2}}.\]
By the uniform control above, we can bound the third term of \eqref{331} by
\[c\sum_{x\in D_{N}^{\varepsilon}}\sum_{k=1}^{n(x)-1}e^{-\widetilde{c}M^{2}(n(x)-k)^{1/2}}\mathbb{P}(h^{D_{N}}(x)-a_{N}\in[b,b'))\leq ce^{-\widetilde{c}M^{2}}\mathbb{E}|\Gamma_{N}^{D}(b)|.\]
By Lemma \ref{3.1}, we complete the proof.
\end{proof}
\par
Now we turn to the truncated second moment estimation.
\par
\begin{proof}[Proof of Proposition \ref{4.2}.]
Writing
\begin{equation}\label{424}
\begin{split}
\mathbb{E}[|\widehat{\Gamma}_{N}^{D,M}(b,b')\cap D_{N}^{\varepsilon}|^{2}]&=\sum_{x,y\in D_{N}^{\varepsilon}}\mathbb{P}(x,y\in\widehat{\Gamma}_{N}^{D,M}(b,b'))\\
&=\sum_{\substack{x,y\in D_{N}^{\varepsilon}\\ \|x-y\|\leq K_{N}^{1/8}}}\mathbb{P}(x,y\in\widehat{\Gamma}_{N}^{D,M}(b,b'))+\sum_{\substack{x,y\in D_{N}^{\varepsilon}\\ \|x-y\|>K_{N}^{1/8}}}\mathbb{P}(x,y\in\widehat{\Gamma}_{N}^{D,M}(b,b')),
\end{split}
\end{equation}
we will separately estimate the two terms on the right-hand side of \eqref{424}. The first term on the right-hand side of \eqref{424} can be bounded by $K_{N}^{1/2}\mathbb{E}|\widehat{\Gamma}_{N}^{D,M}(b,b')|=O(K_{N}^{3/2})$ thanks to Proposition \ref{4.1}.
\par
In order to estimate the second term on the right-hand side of \eqref{424}, for any $x,y$ such that $\|x-y\|>K_{N}^{1/8}$, we need to estimate the probability that $x,y\in\widehat{\Gamma}_{N}^{D,M}(b,b')$. Denote $k:=k(x,y)=(\lceil\log^{+}\|x-y\|\rceil+1)\wedge n(x)$, and let $l=l(x,y)\geq1$ be the minimal such that
\[\Delta^{k-l}(x)\cap\Delta^{k-l}(y)=\varnothing,\quad\Delta^{k-l+1}(x)\cup\Delta^{k-l+1}(y)\subseteq\Delta^{k}(x).\]
Since $n(x)\leq\log N+c(D)$ and $n(y)\geq\log N-c'(\varepsilon)$, we have $l\leq\widetilde{c}(D,\varepsilon)$. Writing
\begin{equation}\label{421}h^{D_{N}}(x)=S_{k}(x)+(S_{k-l}(x)-S_{k}(x))+(S_{0}(x)-S_{k-l}(x))\end{equation}
and
\begin{equation}\label{422}h^{D_{N}}(y)=S_{k}(x)+(S_{k-l}(y)-S_{k}(x))+(S_{0}(y)-S_{k-l}(y)),\end{equation}
Lemma \ref{4.7} implies that the three terms on the right-hand side of \eqref{421} are independent, and the third term is independent of the first two terms on the right-hand side of \eqref{422}. Given $t\in[-M(n(x)-k)^{3/4},M(n(x)-k)^{3/4}]$, $s_{1},s_{2}\in[b,b')$ and $u_{1},u_{2}\in[-n(x)^{3/4},n(x)^{3/4}]$, conditioned on $S_{k}(x)-a_{N}\frac{n(x)-k}{n(x)}=t,S_{k-l}(x)-S_{k}(x)=u_{1},S_{k-l}(y)-S_{k}(x)=u_{2}$, we have
\begin{equation}\label{423}
\begin{split}
&\quad\ \mathbb{P}^{*}(h^{D_{N}}(x)-a_{N}\in\mathrm{d}s_{1},h^{D_{N}}(y)-a_{N}\in\mathrm{d}s_{2})\\
&=\mathbb{P}^{*}(S_{0}(x)-S_{k-l}(x)-a_{N}\tfrac{k}{n(x)}+t+u_{1}\in\mathrm{d}s_{1})\mathbb{P}(S_{0}(y)-S_{k-l}(y)-a_{N}\tfrac{k}{n(x)}+t+u_{2}\in\mathrm{d}s_{2})\\
&\leq\frac{c}{k}\exp\left(-\frac{(\frac{ka_{N}}{n(x)}-t-u_{1}+s_{1})^{2}+(\frac{ka_{N}}{n(x)}-t-u_{2}+s_{2})^{2}}{2k\gamma}\right)\mathrm{d}s_{1}\mathrm{d}s_{2}\\
&\leq\frac{c}{k}\exp\left(-\frac{ka_{N}^{2}}{\gamma n^{2}(x)}+\frac{a_{N}}{\gamma n(x)}(2t+u_{1}+u_{2}-s_{1}-s_{2})\right)\mathrm{d}s_{1}\mathrm{d}s_{2},
\end{split}
\end{equation}
where $\mathbb{P}^{*}$ denotes the conditional probability. In the first inequality of \eqref{423}, we use Lemma \ref{4.6} to replace the variance of $S_{0}(x)-S_{k-l}(x)$ and $S_{0}(y)-S_{k-l}(y)$ by $k\gamma$, which only contributes a harmless multiplicative factor. Next, we will integrate the conditional probability in \eqref{423} with respect to $S_{k-l}(y)-S_{k}(x)$ and $S_{k-l}(x)-S_{k}(x)$ given $S_{k}(x)$. By Lemma \ref{4.6}, there exists $c>0$ and $c''>c'>0$ such that $\mathbb{E}[(S_{k-l}(y)-S_{k}(x))S_{k}(x)]\leq c$ and $\mathrm{var}[S_{k-l}(y)-S_{k}(x)]\in[c',c'']$. For any $t$ with $|t|\leq M(n(x)-k)^{3/4}$, we have
\[\left|\mathbb{E}\left[S_{k-l}(y)-S_{k}(x)\Big|S_{k}(x)-\tfrac{n(x)-k}{n(x)}=t\right]\right|\leq\frac{ct\vee1}{n(x)-k+1}\leq cM,\]
and
\[\mathrm{var}[S_{k-l}(y)-S_{k}(x)|S_{k}(x)]\leq\mathrm{var}[S_{k-l}(y)-S_{k}(x)]\leq c''.\]
For $S_{k-l}(x)-S_{k}(x)$, the corresponding conditional expectation and variance have similar bounds. Since $a_{N}/n(x)$ is bounded, we can get
\[\mathbb{E}\left[e^{\frac{a_{N}}{\gamma n(x)}(S_{k-l}(y)-S_{k}(x))+\frac{a_{N}}{\gamma n(x)}(S_{k-l}(x)-S_{k}(x))}\Big| S_{k}(x)-\tfrac{n(x)-k}{n(x)}=t\right]\leq\widetilde{c}\]
by Cauchy-Schwarz inequality. Applying \eqref{423}, we have
\begin{equation}\label{340}
\begin{split}
&\quad\ \mathbb{P}\left(h^{D_{N}}(x)-a_{N}\in[b,b'),h^{D_{N}}(y)-a_{N}\in[b,b')\Big|S_{k}(x)-\tfrac{n(x)-k}{n(x)}=t\right)\\
&\leq\mathbb{P}\left(|S_{k-l}(x)-S_{k}(x)|\vee|S_{k-l}(y)-S_{k}(x)|>n(x)^{3/4}\Big|S_{k}(x)-\tfrac{n(x)-k}{n(x)}=t\right)+\\
&\quad\ \frac{c}{k}\iint_{[b,b')\times[b,b')}e^{-\frac{ka_{N}^{2}}{\gamma n^{2}(x)}+\frac{a_{N}}{\gamma n(x)}(2t-s_{1}-s_{2})}\mathrm{d}s_{1}\mathrm{d}s_{2}\\
&\leq\frac{c'}{k}\exp\left(-\frac{ka_{N}^{2}}{\gamma n^{2}(x)}+\frac{2a_{N}t}{\gamma n(x)}\right)
\end{split}
\end{equation}
with $c'=c'(b,b',M)$. Note that if $k=n(x)$ then $S_{k}(x)=t=0$ and the condition vanishes. Therefore in this case, the right-hand side of \eqref{340} also bounds the unconditional probability. When $k<n(x)$, we integrate the left-hand side of \eqref{340} with respect to $S_{k}(x)-\tfrac{n(x)-k}{n(x)}$ to get
\begin{equation}\label{345}
\begin{split}
&\quad\ \mathbb{P}(x,y\in\widehat{\Gamma}_{N}^{D,M}(b,b'))\\
&\leq\mathbb{P}\left(h^{D_{N}}(x)-a_{N}\in[b,b'),h^{D_{N}}(y)-a_{N}\in[b,b'),\big|S_{k}(x)-\tfrac{n(x)-k}{n(x)}\big|\leq M(n(x)-k)^{3/4}\right)\\
&\leq\frac{c}{k\sqrt{n(x)-k}}e^{-\frac{ka_{N}^{2}}{\gamma n^{2}(x)}}\int_{|t|\leq M(n(x)-k)^{3/4}}e^{\frac{2a_{N}t}{\gamma n(x)}-\frac{\left(a_{N}\frac{n(x)-k}{n(x)}+t\right)^{2}}{2\gamma(n(x)-k)}}\mathrm{d}t\\
&\leq\frac{c}{k\sqrt{n(x)-k}}e^{-\frac{(k+n(x))a_{N}^{2}}{2\gamma n^{2}(x)}+\widetilde{c}M(n(x)-k)^{3/4}}\leq c''\frac{K_{N}}{N^{4}}\frac{\sqrt{n(x)}}{k\sqrt{n(x)-k}}e^{-\frac{ka_{N}^{2}}{2\gamma n^{2}(x)}+\widetilde{c}M(n(x)-k)^{3/4}}.
\end{split}
\end{equation}
Finally, for the estimation of the second term on the right-hand side of \eqref{424}, we will write:
\[\sum_{\substack{x,y\in D_{N}^{\varepsilon}\\ \|x-y\|>K_{N}^{1/8}}}\mathbb{P}(x,y\in\widehat{\Gamma}_{N}^{D,M}(b,b'))=\sum_{k=k_{N}}^{n}\sum_{\substack{x,y\in D_{N}^{\varepsilon}\\ \|x-y\|>K_{N}^{1/8}\\k(x,y)=k}}\mathbb{P}(x,y\in\widehat{\Gamma}_{N}^{D,M}(b,b'))\]
where we set $n=\max_{x\in D_{N}^{\varepsilon}}n(x)$. Note that, for a given $k$, the sum contains $O(N^{4}e^{4k})$ pairs of $(x,y)$. Since a change in $n$ by an additive constant only contributes a harmless multiplicative factor, by \eqref{345} we have
\begin{equation}\label{425}\sum_{\substack{x,y\in D_{N}^{\varepsilon}\\ \|x-y\|>K_{N}^{1/8}}}\mathbb{P}(x,y\in\widehat{\Gamma}_{N}^{D,M}(b,b'))\leq c'\frac{K_{N}}{N^{4}}\sum_{k=k_{N}}^{n}\frac{\sqrt{n}}{k\sqrt{n-k+1}}N^{4}e^{4k-\frac{ka_{N}^{2}}{2\gamma n^{2}}+\widetilde{c}M(n-k)^{3/4}}.\end{equation}
Since $\frac{a_{N}^{2}}{2\gamma n^{2}}\to4\lambda^{2}<4$ as $N\to\infty$, the exponent on the right-hand side of \eqref{425} grows linearly with $k$. Then we can dominate the sum by the $k=n$ term. As $n=\log N+O(1)$, the right-hand side is $O(K_{N}^{2})$ by some simple calculation. In order to complete the proof, we only need to rescale the domain and show the diameter dependence since all bounds above are regardless of the diameter of the domain.
\par
We now rescale the domain with a factor $r<1$. Assuming $D':=r^{-1}D$, note that the lattice approximation and the scale sequence changes as $D'_{N}:=D_{[N/r]-j},a'_{N}:=a_{[N/r]-j},k'_{N}:=k_{[N/r]-j}$ for some $j$, where $j$ is appropriately chosen in $\{0,1,\cdots,[r^{-1}]\}$. We then get
\[\widehat{\Gamma}_{[N/r]-j}^{D,M}(b)=\widehat{\Gamma}_{N}^{D',M}(b).\]
Note that
\[K'_{N}:=\frac{N^{4}}{\sqrt{\log N}}e^{-\frac{(a'_{N})^{2}}{2\gamma\log N}}=(r^{4+4\lambda^{2}}+o(1))K_{[N/r]-j},\]
as all integers can be represented by $[N/r]-j$ with some $N$ and $j$, the claim for $D$ follows from the claim of $D'$, which completes the proof.
\end{proof}
\section{Existence and factorization of subsequential limits}
In this part, we will prove that subsequential weak limit of the measures $\{\eta_{N}^{D}\}$ exists, and can be factorized as in \eqref{fac}. Recall that $\widehat{\eta}_{N}^{D,M}$ is the truncated point measure defined in \eqref{tpm}. Since $\overline{D}\times(\mathbb{R}\cup\{\infty\})$ is separable, we have:
\par
\begin{proposition}\label{3.4}
For $\lambda\in(0,1)$ and large $M$, $\{\widehat{\eta}_{N}^{D,M}:N\geq1\}$ is tight with respect to the topology of vague convergence on the space of Radon measures on $\overline{D}\times(\mathbb{R}\cup\{\infty\})$. Moreover, for any subsequential weak limit $\widehat{\eta}^{D,M}$ of these measures and any $b\in\mathbb{R}$,
\begin{equation}\label{411}\mathbb{P}(\widehat{\eta}^{D,M}(D\times[b,\infty))<\infty)=1,\end{equation}
and, for any non-empty open set $A\subset D$,
\begin{equation}\label{412}\mathbb{P}(\widehat{\eta}^{D,M}(A\times[b,\infty))>0)>0.\end{equation}
Moreover, we have $\widehat{\eta}^{D,M}(A\times\mathbb{R})=0$ a.s. for any $A$ with $\mu(A)=0$, with $\mu$ the Lebesgue measure.
\end{proposition}
\begin{proof}
By definition, $\widehat{\eta}_{N}^{D,M}(\overline{D}\times[b,\infty))=|\widehat{\Gamma}_{N}^{D,M}(b)|\leq|\Gamma_{N}^{D}(b)|$. Then by Lemma \ref{3.1}, the family of measures $\{\widehat{\eta}_{N}^{D,M}(\overline{D}\times[b,\infty)):N\geq1\}$ is tight for each $b\in\mathbb{R}$. Thus for any function $f:\overline{D}\times(\mathbb{R}\cup\{\infty\})\to\mathbb{R}$ which is continuous with compact support, $\{\langle\widehat{\eta}_{N}^{D,M},f\rangle:N\geq1\}$ is tight, which completes the proof of the first statement.
\par
Now we let $\widehat{\eta}^{D,M}$ be a subsequential weak limit of the measures $\{\widehat{\eta}_{N}^{D,M}:N\geq1\}$. By Fatou's lemma, Lemma \ref{3.1} and the domination of $\widehat{\eta}_{N}^{D,M}$ by $\eta_{N}^{D}$, we have $\mathbb{E}\widehat{\eta}^{D,M}(\overline{D}\times[b,\infty))<\infty$ for each $b\in\mathbb{R}$. Lemma \ref{3.1} also gives $\widehat{\eta}^{D,M}(A\times\mathbb{R})=0$ a.s. for any $A$ with $\mu(A)=0$.
\par
For \eqref{412}, denote $X_{N}:=\widehat{\eta}_{N}^{D,M}(A\times[b,\infty))$. Note that $\{X_{N}:N\geq1\}$ is uniformly integrable since $\sup_{N\geq1}\mathbb{E}X_{N}^{2}<\infty$ thanks to Proposition \ref{4.2}. Lemma \ref{3.2} implies that $\inf_{N\geq1}\mathbb{E}X_{N}^{2}>0$. Obviously $\exp(\frac{4\lambda^{2}}{\gamma}s_{D}(x))>0$, then any distributional limit of $X_{N}$ has positive expectation for any open set $A$ such that $\varnothing\neq A\subset D$, which completes the proof.
\end{proof}
\par
The claim \eqref{412} shows that every subsequential weak limit $\widehat{\eta}^{D,M}$ of measures $\{\widehat{\eta}_{N}^{D,M}:N\geq1\}$ has positive total mass with positive probability when $M$ is large enough. In light of the domination of $\widehat{\eta}_{N}^{D,M}$ by $\eta_{N}^{D}$, \eqref{411} and \eqref{412} remain valid with $\widehat{\eta}_{N}^{D,M}$ replaced by $\eta^{D}$.
\par
Then, for any $b\in\mathbb{R}$ and function $f:\overline{D}\times(\mathbb{R}\cup\{\infty\})\to\mathbb{R}$, define
\[f_{b}(x,h):=f(x,h+b)e^{-\pi\lambda b}.\]
The following proposition will lead to the factorization property:
\par
\begin{proposition}\label{3.5}
Let $\lambda\in(0,1)$, and $\eta^{D}$ be any subsequential limit of $\{\eta_{N}^{D}:N\geq1\}$. Then for any $b\in\mathbb{R}$ and any $f:\overline{D}\times(\mathbb{R}\cup\{\infty\})\to\mathbb{R}$ of the form $f(x,h)=1_{A}(x)1_{[0,\infty)}(h)$ with $A\subset D$ open,
\[\langle\eta^{D},f_{b}\rangle=\langle\eta^{D},f\rangle,\quad\mathrm{a.s.}\]
\end{proposition}
\par
In order to prove the proposition, we need the following lemma:
\par
\begin{lemma}\label{4.3}
Let $\lambda\in(0,1)$. For each open $A\subset D$ and each $b\in\mathbb{R}$, denoting $A_{N}:=\{x\in\mathbb{Z}^{4}:x/N\in A\}$, we have
\[\lim_{N\to\infty}\frac{1}{K_{N}}\mathbb{E}\left||\widehat{\Gamma}_{N}^{D,M}(0)\cap A_{N}|-e^{\pi\lambda b}|\widehat{\Gamma}_{N}^{D,M}(b)\cap A_{N}|\right|=0.\]
\end{lemma}
\begin{proof}
By Cauchy-Schwarz inequality, it suffices to prove
\begin{equation}\label{431}\lim_{N\to\infty}\frac{1}{K_{N}^{2}}\mathbb{E}\left(|\widehat{\Gamma}_{N}^{D,M}(0)\cap A_{N}|-e^{\pi\lambda b}|\widehat{\Gamma}_{N}^{D,M}(b)\cap A_{N}|\right)^{2}=0.\end{equation}
The expectation in \eqref{431} can be written as
\begin{equation}\label{432}
\begin{split}
\sum_{x,y\in A_{N}}\mathbb{E}&[(1_{\{h^{D_{N}}(x)\geq a_{N}\}}-e^{\pi\lambda b}1_{\{h^{D_{N}}(x)\geq a_{N}+b\}})\\
&(1_{\{h^{D_{N}}(y)\geq a_{N}\}}-e^{\pi\lambda b}1_{\{h^{D_{N}}(y)\geq a_{N}+b\}})1_{T_{N,M}(x)}1_{T_{N,M}(y)}].
\end{split}
\end{equation}
We first consider the terms corresponding to $x,y$ such that
\begin{equation}\label{430}m:=[\log\|x-y\|]\end{equation}
satisfies
\begin{equation}\label{433}m\geq\frac{3}{2}k_{N},\quad\|x-y\|\in[e^{m}+2e^{k_{N}},e^{m+1}-2e^{k_{N}}].\end{equation}
Let $\mathscr{F}:=\sigma(h^{D_{N}}(z):z\in D_{N}\backslash(\Delta^{k_{N}}(x)\cup\Delta^{k_{N}}(y)))$. \eqref{430} and \eqref{433} ensure that
\[\Delta^{k_{N}}(x)\cap\partial_{2}\Delta^{k}(y)=\varnothing,\quad k=k_{N},\ldots,n(y),\]
and
\[\Delta^{k_{N}}(y)\cap\partial_{2}\Delta^{k}(x)=\varnothing,\quad k=k_{N},\ldots,n(x),\]
which implies that $T_{N,M}(x)$ and $T_{N,M}(y)$ are both $\mathscr{F}$-measurable, since $S_{k}(x)$ (resp.\ $S_{k}(y)$) only depends on $\{h^{D_{N}}(z):z\in\partial_{2}\Delta^{k}(x)\}$ (resp.\ $\{h^{D_{N}}(z):z\in\partial_{2}\Delta^{k}(y)\}$).
\par
Therefore, the term in \eqref{432} corresponding to such $x,y$ can be written as
\begin{equation}\label{439}
\begin{split}
&\mathbb{E}[(\mathbb{P}(h^{D_{N}}(x)\geq a_{N}|\mathscr{F})-e^{\pi\lambda b}\mathbb{P}(h^{D_{N}}(x)\geq a_{N}+b|\mathscr{F}))1_{T_{N,M}(x)}\\
&(\mathbb{P}(h^{D_{N}}(y)\geq a_{N}|\mathscr{F})-e^{\pi\lambda b}\mathbb{P}(h^{D_{N}}(y)\geq a_{N}+b|\mathscr{F}))1_{T_{N,M}(y)}].
\end{split}
\end{equation}
Thanks to Theorem \ref{gmp}, we can write $h^{D_{N}}(x)=(S_{0}(x)-S_{k_{N}}(x))+S_{k_{N}}(x)$, where the two terms on the right-hand side are independent of each other, and $(S_{0}(x)-S_{k_{N}}(x))$ has the law of the membrane model in $\Lambda_{e^{k_{N}}}(x)$. Using a similar decomposition for $h^{D_{N}}(y)$, and thanks to the fact that $S_{0}(x)-S_{k_{N}}(x)$ and $S_{0}(y)-S_{k_{N}}(y)$ are both $\mathscr{F}$-measurable, then \eqref{439} is bounded by
\[\mathbb{E}\left[1_{T_{N,M}(x)}F\left(S_{k_{N}}(x)-a_{N}\tfrac{\log N-k_{N}}{\log N}\right)1_{T_{N,M}(y)}F\left(S_{k_{N}}(y)-a_{N}\tfrac{\log N-k_{N}}{\log N}\right)\right],\]
where
\[F(u)=\left|\mathbb{P}\left(\widetilde{h}(0)\geq a_{N}\tfrac{k_{N}}{\log N}-u\right)-e^{\pi\lambda b}\mathbb{P}\left(\widetilde{h}(0)\geq a_{N}\tfrac{k_{N}}{\log N}-u+b\right)\right|\]
with $\widetilde{h}$ denoting the membrane model on $\Lambda_{e^{k_{N}}}$. Since $a_{N}k_{N}/\log N\sim2\lambda\sqrt{2\gamma}k_{N}$ and $\mathbb{E}[\widetilde{h}^{2}(0)]=\gamma k_{N}+O(1)$, if we assume that $|u|<k_{N}^{7/8}$, thanks to \eqref{321}, we then have
\[\frac{\mathbb{P}\left(\widetilde{h}(0)\geq a_{N}\tfrac{k_{N}}{\log N}-u\right)}{\mathbb{P}\left(\widetilde{h}(0)\geq a_{N}\tfrac{k_{N}}{\log N}-u+b\right)}\sim\frac{a_{N}\tfrac{k_{N}}{\log N}-u+b}{a_{N}\tfrac{k_{N}}{\log N}-u}\exp\left(\frac{-2b(a_{N}\tfrac{k_{N}}{\log N}-u)+b^{2}}{2\mathbb{E}[\widetilde{h}^{2}(0)]}\right)\sim e^{\pi\lambda b}\]
as $N\to\infty$. Since on $T_{N,M}(x)$ we have $|S_{k_{N}}(x)-a_{N}\tfrac{\log N-k_{N}}{\log N}|<k_{N}^{7/8}$ when $N$ is sufficiently large, we have $F(u)=o(1)\mathbb{P}\left(\widetilde{h}(0)\geq a_{N}\tfrac{k_{N}}{\log N}-u\right)$ with $o(1)\to0$ as $N\to\infty$ uniformly in $|u|<k_{N}^{7/8}$. Similarly, on $T_{N,M}(y)$ we have $|S_{k_{N}}(y)-a_{N}\tfrac{\log N-k_{N}}{\log N}|<k_{N}^{7/8}$ when $N$ is sufficiently large. Then we can show that the expectation is bounded by
\[o(1)\mathbb{P}(h^{D_{N}}(x)\geq a_{N},h^{D_{N}}(y)\geq a_{N},T_{N,M}(x),T_{N,M}(y)).\]
As in Proposition \ref{4.2} the sum of the corresponding terms is thus at most $o(K_{N}^{2})$ in this case.
\par
We now consider the remaining terms in the sum. Recall \eqref{430}, if $x,y$ satisfy $m<\frac{3}{2}k_{N}$, then we can bound the corresponding term by $4e^{2\pi\lambda(b\vee0)}\mathbb{P}(h^{D_{N}}(x)\geq a_{N}+b\wedge0)$. Thanks to Lemma \ref{3.1}, the sum of these terms is at most $o(K_{N}^{2})$. If $x,y$ satisfy $m\geq\frac{3}{2}K_{N}$ but not the second condition in \eqref{433}, then we can bound the corresponding term by $4e^{2\pi\lambda(b\vee0)}\mathbb{P}(h^{D_{N}}(x)\geq a_{N}-b\wedge0,h^{D_{N}}(y)\geq a_{N}-b\wedge0)$. Note that the number of such pairs for a given $m$ is at most of order $N^{4}e^{3m+k_{N}}=o(1)N^{4}e^{2m}$, repeating the argument as in \eqref{425}, the sum in this case is at most $o(K_{N}^{2})$ as well, and thus completing the proof.
\end{proof}
\par
We can now prove Proposition \ref{3.5}.
\par
\begin{proof}[Proof of Proposition \ref{3.5}] Let $f(x,h)=1_{A}(x)1_{[0,\infty)}(h)$. Lemma \ref{4.3} can be rephrased as
\[\lim_{N\to\infty}\mathbb{E}\left|\langle\widehat{\eta}_{N}^{D,M},f\rangle-\langle\widehat{\eta}_{N}^{D,M},f_{b}\rangle\right|=0,\quad b\in\mathbb{R}.\]
Taking the distributional limit of $\langle\widehat{\eta}_{N}^{D,M},(f-f_{b})\rangle$, by Fatou's lemma, we have $\langle\widehat{\eta}^{D,M},f-f_{b}\rangle=0$ a.s. Then letting $M\to\infty$, we complete the proof of the proposition.
\end{proof}
\par
The proposition now leads to the following factorization property:
\begin{proposition}\label{3.7}
Let $\eta^{D}$ be any subsequential limit of $\eta_{N}^{D}$ along a deterministic subsequence $\{N_{k}\}$. Then, with probability one,
\begin{equation}\label{fac}\eta^{D}(\mathrm{d}x\mathrm{d}h)=Z_{\lambda}^{D}(\mathrm{d}x)\otimes e^{-\pi\lambda h}\mathrm{d}h\end{equation}
with $Z_{\lambda}^{D}$ a finite random Borel measure on $\overline{D}$.
\end{proposition}
\par
\begin{proof}
Thanks to Proposition \ref{3.5}, we have $\langle\eta^{D},f-f_{b}\rangle=0$ a.s. for each function $f:\overline{D}\times\mathbb{R}\to\mathbb{R}$ of the form $f(x,h)=1_{A}(x)1_{[0,\infty)}(h)$ with $A\subset D$ open and each $b\in\mathbb{R}$. For $A\subset\overline{D}$ Borel, define $Z_{\lambda}^{D}(A):=\pi\lambda\eta^{D}(A\times[0,\infty))$, which is a finite random Borel measure on $D$. If $A$ is open then we have
\begin{align*}
\eta^{D}(A\times[b,\infty))=e^{-\pi\lambda b}\langle\eta^{D},f_{b}\rangle&=e^{-\pi\lambda b}\langle\eta^{D},f\rangle\\
&=(\pi\lambda)^{-1}e^{-\pi\lambda b}Z_{\lambda}^{D}(A)=\int_{A\times[b,\infty)}Z_{\lambda}^{D}(\mathrm{d}x)e^{-\pi\lambda h}\mathrm{d}h
\end{align*}
holds almost surely, where the null set may depend on $A$ or $b$. In order to complete the proof, we can choose a common null set in the $\pi$-system $\{A\times[b,\infty):A\subset D$ open dyadic cube, $b\in\mathbb{Q}\}$ since it is countable. Note that the product Borel $\sigma-$algebra on $D\times\mathbb{R}$ can be generated from the $\pi-$system above, and $\eta^{D}(\partial D\times\mathbb{R})=Z_{\lambda}^{D}(\partial D)=0$ thanks to Proposition \ref{3.4}, the equality then extends to all the Borel sets on $\overline{D}\times\mathbb{R}$.
\end{proof}
\section{Uniqueness of the subsequential limit}
By the previous sections, any subsequential limit $\eta^{D}$ of $\{\eta_{N}^{D}:N\geq1\}$ can be factorized as in Proposition \ref{3.7}. In this section, our aim is to show the uniqueness of the measure $Z_{\lambda}^{D}$ on the right-hand side of \eqref{fac}, and thus also the subsequential limit $\eta^{D}$, which particularly shows that $\eta_{N}^{D}$ converges to the same distributional limit.
\par
Recall $\widetilde{G}^{D}(x,y)=\widetilde{G}_{y}^{D}(x)$ as in \eqref{cgf}, and $s_{D}(x)$ as in \eqref{sdx}. Then, for two domains $\widetilde{D}\subset D$, the function
\[C^{D,\widetilde{D}}(x,y):=\widetilde{G}^{D}(x,y)-\widetilde{G}^{\widetilde{D}}(x,y)\]
defines a symmetric, positive semi-definite function from $\widetilde{D}^{2}$ to $\mathbb{R}$. This allow us to define $\{\Phi^{D,\widetilde{D}}(x):x\in\widetilde{D}\}$ to be a centered Gaussian field with covariance $C^{D,\widetilde{D}}$, which has smooth sample paths a.s. The following proposition shows that $Z_{\lambda}^{D}$ satisfies some properties (in particular the fifth one, the Gibbs-Markov property in the limit) which, as we are going to show in Proposition \ref{3.11}, are sufficient to determine its law uniquely.
\begin{proposition}\label{3.8}
Let $\lambda\in(0,1)$ and let $\{\eta^{D}:D\in\mathfrak{D}\}$ be the subsequential limit of $\{\eta_{N}^{D}:N\geq1\}$ along a deterministic subsequence $\{N_{k}\}$ for $D\in\mathfrak{D}$. Let $Z_{\lambda}^{D}$ be the measure factorized from $\eta^{D}$, then $\{Z_{\lambda}^{D}:D\in\mathfrak{D}\}$ satisfies:
\par
(1) $Z_{\lambda}^{D}$ is supported on $D$; i.e., $Z_{\lambda}^{D}(\mathbb{R}^{4}\backslash D)=0$ a.s.
\par
(2) If $A\subset D$ is measurable with $\mu(A)=0$, then $Z_{\lambda}^{D}(A)=0$ a.s.
\par
(3) For each measurable $A\subset D$,
\[\mathbb{E}Z_{\lambda}^{D}(A)=\frac{\sqrt{\pi}}{4}\int_{A}e^{\frac{4\lambda^{2}}{\gamma}s_{D}(x)}\mathrm{d}x.\]
(4) If $D,\widetilde{D}\in\mathfrak{D}$ such that $D\cap\widetilde{D}=\varnothing$, then
\[Z_{\lambda}^{D\cup\widetilde{D}}(\mathrm{d}x)\overset{\mathrm{law}}{=}Z_{\lambda}^{D}(\mathrm{d}x)+Z_{\lambda}^{\widetilde{D}}(\mathrm{d}x),\]
where $Z_{\lambda}^{D}$ and $Z_{\lambda}^{\widetilde{D}}$ are independent.
\par
(5) If $D,\widetilde{D}\in\mathfrak{D}$ such that $\widetilde{D}\subset D$ and $\mu(D\backslash\widetilde{D})=0$, then
\[Z_{\lambda}^{D}(\mathrm{d}x)\overset{\mathrm{law}}{=}e^{\pi\lambda\Phi^{D,\widetilde{D}}(x)}Z_{\lambda}^{\widetilde{D}}(\mathrm{d}x),\]
where $\Phi^{D,\widetilde{D}}$ is independent of $Z_{\lambda}^{\widetilde{D}}$.
\par
(6) $Z_{\lambda}^{a+D}(a+\mathrm{d}x)\overset{\mathrm{law}}{=}Z_{\lambda}^{D}(\mathrm{d}x)$ for each $a\in\mathbb{R}^{4}$.
\end{proposition}
\begin{proof}
Property (1) is obvious by the definition of the point measure. Property (2) is a direct consequence of Lemma \ref{3.2}. Property (3) holds for all $A\subset D$ open since
\[\mathbb{E}Z_{\lambda}^{D}(A)=\pi\lambda\mathbb{E}\langle\eta^{D},f\rangle=\pi\lambda\lim_{N\to\infty}\mathbb{E}\langle\eta_{N}^{D},f\rangle=\frac{\sqrt{\pi}}{4}\int_{A}e^{\frac{4\lambda^{2}}{\gamma}s_{D}(x)}\mathrm{d}x\]
with $f(x,h)=1_{A}(x)1_{[0,\infty)}(h)$ thanks to Lemma \ref{3.2}. Property (4) is trivial by decompositing $h^{D_{N}\cup\widetilde{D}_{N}}$ by the independent sum of $h^{D_{N}}$ and $h^{\widetilde{D}_{N}}$. Property (6) is trivially true by translation invariance.
\par
Now we turn to property (5). Let $D,\widetilde{D}\in\mathfrak{D}$ with $\mu(D\backslash\widetilde{D})=0$. Writing $h^{D_{N}}\overset{\mathrm{law}}{=}h^{\widetilde{D}_{N}}+\varphi^{D_{N},\widetilde{D}_{N}}$ thanks to Theorem \ref{gmp}, if $f:\overline{D}\times\mathbb{R}\to\mathbb{R}$ is continuous and compactly-supported in $\widetilde{D}$, then
\[\langle\eta_{N}^{D},f\rangle\overset{\mathrm{law}}{=}\langle\eta_{N}^{\widetilde{D}},f_{\varphi}\rangle\]
where
\[f_{\varphi}(x,h)=f(x,h+\varphi^{D_{N},\widetilde{D}_{N}}([xN]))\]
with $\varphi^{D_{N},\widetilde{D}_{N}}$ independent of $h^{\widetilde{D}_{N}}$. Using Theorem \ref{gfa}, we have:
\[G^{D_{N}}([xN],[yN])-G^{\widetilde{D}_{N}}([xN],[yN])\to C^{D,\widetilde{D}}(x,y)\]
as $N\to\infty$ locally uniformly in $x,y\in\widetilde{D}$. Note that $\varphi^{D_{N},\widetilde{D}_{N}}$ is a centered Gaussian process with covariance kernel $G^{D_{N}}-G^{\widetilde{D}_{N}}$ and $\Phi^{D,\widetilde{D}}$ is a centered Gaussian process with covariance kernel $C^{D,\widetilde{D}}(x,y)$. By the same line of the proof in \cite[Theorem B.14]{ref13}, for any $\delta>0$, there exists a coupling of $\varphi^{D_{N},\widetilde{D}_{N}}$ with $\Phi^{D,\widetilde{D}}$ such that
\[\lim_{N\to\infty}\mathbb{P}\left(\sup_{x\in\widetilde{D},d(x,D^{c})>\delta}|\varphi^{D_{N},\widetilde{D}_{N}}([xN])-\Phi^{D,\widetilde{D}}(x)|>\delta\right)=0.\]
Since $f$ is continuous with compact support, letting $f^{\Phi}(x,h)=f(x,h+\Phi^{D,\widetilde{D}}(x))$ with $\Phi^{D,\widetilde{D}}$ independent of $\eta_{N}^{D}$, we thus have
\[\langle\eta_{N}^{\widetilde{D}},f_{\varphi}\rangle\overset{\mathrm{law}}{=}\langle\eta_{N}^{\widetilde{D}},f_{\Phi}\rangle+o(1)\]
with $o(1)\to0$ in probability as $N\to\infty$. Since $x\to\Phi^{D,\widetilde{D}}(x)$ is continuous on $\widetilde{D}$ almost surely, for any subsequential limits $\eta^{D}$ (resp. $\eta^{\widetilde{D}}$) of $\{\eta_{N}^{D}:N\geq1\}$ (resp. $\{\eta_{N}^{\widetilde{D}}:N\geq1\}$), we thus get
\[\langle\eta^{D},f\rangle\overset{\mathrm{law}}{=}\langle\eta^{\widetilde{D}},f_{\Phi}\rangle.\]
Then by Proposition \ref{3.7}, we have
\[\langle\eta^{\widetilde{D}},f_{\Phi}\rangle=\int_{D\times\mathbb{R}}f(x,h+\Phi^{D,\widetilde{D}}(x))Z_{\lambda}^{\widetilde{D}}(\mathrm{d}x)e^{-\pi\lambda h}\mathrm{d}h=\int_{D\times\mathbb{R}}f(x,h)Z_{\lambda}^{\widetilde{D}}(\mathrm{d}x)e^{-\pi\lambda(h-\Phi^{D,\widetilde{D}}(x))}\mathrm{d}h.\]
As the equality holds for any continuous $f:D\times\mathbb{R}\to\mathbb{R}$ supported on $\widetilde{D}$, we have
\[Z_{\lambda}^{\widetilde{D}}(\mathrm{d}x)e^{\pi\lambda\Phi^{D,\widetilde{D}}(x)}\overset{\mathrm{law}}{=}Z_{\lambda}^{D}(\mathrm{d}x),\]
which is the desired result.
\end{proof}
\par
For uniqueness of $Z_{\lambda}^{D}$, we claim that $Z_{\lambda}^{S_{n}}$ has the same law as the limit of the measures
\begin{equation}\label{ym}Y_{m}^{S_{n}}(\mathrm{d}x):=\frac{\sqrt{\pi}}{4}\sum_{j=1}^{16^{m}}e^{\pi\lambda\Phi^{S_{n},\widetilde{S}_{n,m}}(x)+\frac{4\lambda^{2}}{\gamma}s_{S_{n+m,j}}(x)}1_{S_{n+m,j}}(x)\mathrm{d}x\end{equation}
as $m\to\infty$, where $S_{n}=(0,2^{-n})^{4}$, and $\widetilde{S}_{n,m}:=\cup_{j=1}^{16^{m}}S_{n+m,j}$ with $S_{n+m,j}(j=1,\cdots,16^{m})$ the collection of the translations of $S_{n+m}$ decomposing from $S_{n}$. The following lemma gives the (weak) limit of $Y_{m}^{S_{n}}$.
\begin{lemma}\label{3.10}
For each $\lambda\geq0$, there exists an a.s.\ finite random measure $Y_{\infty}^{S_{n}}$, such that for each bounded, measurable $f:S_{n}\to\mathbb{R}$,
\[\lim_{m\to\infty}\langle Y_{m}^{S_{n}},f\rangle=\langle Y_{\infty}^{S_{n}},f\rangle\quad\mathrm{a.s.}\]
\end{lemma}
\begin{proof}
We can write
\[\Phi^{S_{n},\widetilde{S}_{n,m}}=\sum_{j=1}^{m}\Phi^{\widetilde{S}_{n,j-1},\widetilde{S}_{n,j}}(x),\]
where $\widetilde{S}_{n,0}:=S_{n}$ by the covariance structure of the field, and the terms on the right-hand side are independent. Then we immediately get that $\langle Y_{m}^{S_{n}},f\rangle$ is a martingale adapted to the filtration
\[\mathscr{F}_{m}:=\sigma(\Phi^{\widetilde{S}_{n,j-1},\widetilde{S}_{n,j}}(x):x\in S'_{n},j=1,2,\cdots,m)\quad\mathrm{where}\quad S'_{n}:=\bigcap_{m\geq1}\widetilde{S}_{n,m},\] since
\[\sum_{j=1}^{16^{m}}e^{\frac{4\lambda^{2}}{\gamma}s_{S_{n+m,j}}(x)}1_{S_{n+m,j}}(x)=e^{\frac{4\lambda^{2}}{\gamma}s_{\widetilde{S}_{n,m}}(x)},\quad x\in\widetilde{S}_{n,m},\]
and, for any $\widetilde{D}\subset{D}$,
\[\mathbb{E}e^{\pi\lambda\Phi^{D,\widetilde{D}}(x)}=e^{\frac{1}{2}\pi^{2}\lambda^{2}C^{D,\widetilde{D}}(x)}=e^{\frac{4\lambda^{2}}{\gamma}(s_{D}(x)-s_{\widetilde{D}}(x))},\quad x\in\widetilde{D}.\]
By separately considering the positive and negative part of $f$, and applying martingale convergence theorem to both parts, we can get
\[L(f):=\lim_{m\to\infty}\langle Y_{m}^{S_{n}},f\rangle\]
exists almost surely.
\par
It remains to show that the limit is an integral of $f$ with respect to a random measure. We repeat the argument in \cite[Lemma 3.10]{ref1}. We fix a countable dense subset $A\in C(\overline{D})$. By Fatou's lemma, we have
\[\mathbb{E}|L(f)|\leq c\int_{S_{n}}|f(x)|e^{\frac{4\lambda^{2}}{\gamma}s_{D}(x)}\mathrm{d}x,\quad\forall f\in C(\overline{D}).\]
Therefore the linear functional $f\mapsto L(f)$ is almost surely well-defined for all $f\in A$ simultaneously and is bounded on $A$. Then $f\mapsto L(f)$ can be extended uniquely to a linear functional $f\mapsto\overline{L}(f)$ on $C(\overline{D})$ such that $L(f)=\overline{L}(f)$ almost surely for any $f\in C(\overline{D})$. The existence of $Y_{\infty}^{S_{n}}$ then follows from Riesz representation theorem, and thus complete the proof.
\end{proof}
\par
The following Proposition give the uniqueness of the law of $Z_{\lambda}^{D}$:
\begin{proposition}\label{3.11}
For $\lambda\in(0,1)$ and any $n\in\mathbb{Z}$,
\[Z_{\lambda}^{S_{n}}(\mathrm{d}x)\overset{\mathrm{law}}{=}Y_{\infty}^{S_{n}}(\mathrm{d}x).\]
\end{proposition}
\begin{proof}
For any bounded, measurable function $f:S_{n}\to[0,\infty)$, we are going to show that
\begin{equation}\label{541}\mathbb{E}e^{-\langle Z_{\lambda}^{S_{n}},f\rangle}=\mathbb{E}e^{-\langle Y_{\infty}^{S_{n}},f\rangle},\end{equation}
which will be done by separately proving the ``$\geq$'' and ``$\leq$'' inequalities.
\par
\textit{Proof of ``$\geq$''.}\quad Thanks to properties (4) and (5) in Proposition \ref{3.8}, we can write
\[Z_{\lambda}^{S_{n}}(\mathrm{d}x)=\sum_{j=1}^{16^{m}}e^{\pi\lambda\Phi^{S_{n},\widetilde{S}_{n,m}}(x)}1_{S_{n+m,j}}(x)Z_{\lambda}^{S_{n+m,j}}(\mathrm{d}x),\]
where $Z_{\lambda}^{S_{n+m,j}},1\leq j\leq16^{m}$ are independent of others as well as of $\Phi^{S_{n},\widetilde{S}_{n,m}}$. Then we have
\begin{equation}\label{540}\mathbb{E}[\langle Z_{\lambda}^{S_{n}},f\rangle|\sigma(\Phi^{S_{n},\widetilde{S}_{n,m}})]=\langle Y_{m}^{S_{n}},f\rangle.\end{equation}
By Jensen's inequality, we have
\[\mathbb{E}e^{-\langle Z_{\lambda}^{S_{n}},f\rangle}\geq\mathbb{E}e^{-\langle Y_{m}^{S_{n}},f\rangle}.\]
Let $m\to\infty$, by the bounded convergence theorem, we complete the proof of ``$\geq$'' in \eqref{541}.
\par
\textit{Proof of ``$\leq$''.}\quad We need to introduce an additional truncation: For $\delta\in(0,1/2)$, let $S_{k}^{\delta}$ be the translate of $(\delta2^{-k},(1-\delta)2^{-k})^{4}$ that has the same center as $S_{k}$. Let $S_{n+m,j}^{\delta}$ denote the truncated version of $S_{n+m,j}$, then let $\widetilde{S}_{n,m}^{\delta}:=\cup_{j=1}^{16^{m}}S_{n+m,j}^{\delta}$ and define
\[f_{m,\delta}(x):=1_{\widetilde{S}_{n,m}^{\delta}}(x)f(x).\]
The positivity of $f$ then implies $f_{m,\delta}\uparrow f$ as $\delta\downarrow0$. By Markov inequality, for any $\varepsilon>0$ we have $\mathbb{P}(Y_{m}^{S_{n}}(S_{n}\backslash\widetilde{S}_{n,m}^{\delta})\geq\varepsilon)\leq c\delta/\varepsilon$. Note that
\[\langle Y_{m}^{S_{n}},f_{m,\delta}\rangle\leq\langle Y_{m}^{S_{n}},f\rangle\leq\langle Y_{m}^{S_{n}},f_{m,\delta}\rangle+\|f\|Y_{m}^{S_{n}}(S_{n}\backslash\widetilde{S}_{n,m}^{\delta}),\]
we conclude that
\[\lim_{\delta\to0}\limsup_{m\to\infty}\mathbb{E}|\langle Y_{m}^{S_{n}},f_{m,\delta}\rangle-\langle Y_{\infty}^{S_{n}},f\rangle|=0,\]
which implies that $\langle Y_{m}^{S_{n}},f_{m,\delta}\rangle\to\langle Y_{\infty}^{S_{n}},f\rangle$ almost surely as $\delta\to0$ and $m\to\infty$. So we can work with $f_{m,\delta}$ instead of $f$. We now define
\[\widehat{Z}_{\lambda}^{D,M}(A):=\pi\lambda\widehat{\eta}^{D,M}(A\times[0,\infty))\]
as the ``truncated'' version of $Z_{\lambda}^{D}$, and
\[\widetilde{Z}_{\lambda}^{S_{n},M,m}(\mathrm{d}x)=\sum_{j=1}^{16^{m}}e^{\pi\lambda\Phi^{S_{n},\widetilde{S}_{n,m}}(x)}1_{S_{n+m,j}}(x)\widehat{Z}_{\lambda}^{S_{n+m,j},M}(\mathrm{d}x).\]
Thanks to Proposition \ref{4.1} we have
\[\widehat{Z}_{\lambda}^{D,M}(A)\leq Z_{\lambda}^{D}(A),\quad\widehat{Z}_{\lambda}^{D,M}(A)\uparrow Z_{\lambda}^{D}(A)(M\to\infty).\]
For any bounded, measurable $f:S_{n}\to\mathbb{R}_{+}$ and $\delta>0$, we then have
\[\mathbb{E}e^{-\langle Z_{\lambda}^{S_{n}},f\rangle}\leq\mathbb{E}e^{-\langle\widetilde{Z}_{\lambda}^{S_{n},M,m},f\rangle}\leq\mathbb{E}e^{-\langle\widetilde{Z}_{\lambda}^{S_{n},M,m},f_{m,\delta}\rangle}.\]
Now set
\[\widetilde{X}_{j}:=\int_{S_{n+m,j}^{\delta}}e^{\pi\lambda\Phi^{S_{n},\widetilde{S}_{n,m}}(x)}f_{m,\delta}(x)\widehat{Z}_{\lambda}^{S_{n+m,j},M}(\mathrm{d}x).\]
By the fact that $\langle\widetilde{Z}_{\lambda}^{S_{n},M,m},f_{m,\delta}\rangle=\sum_{j=1}^{16^{m}}\widetilde{X}_{j}$, and thanks to Lemma \ref{lem:rev}, we then get for each $\varepsilon>0$,
\begin{equation}\label{542}\mathbb{E}e^{-\langle\widetilde{Z}_{\lambda}^{S_{n},M,m},f_{m,\delta}\rangle}\leq\mathbb{E}\left[\exp\left(-e^{-\varepsilon}\sum_{j=1}^{16^{m}}\mathbb{E}[\widetilde{X}_{j}1_{\{\widetilde{X}_{j}\leq\varepsilon\}}|\Phi^{S_{n},\widetilde{S}_{n,m}}]\right)\right].\end{equation} Then we need the following lemma:
\par
\begin{lemma}\label{3.12}
Suppose $\lambda\in(0,1),\delta\in(0,1/2)$, for each $\varepsilon>0$ we have:
\begin{equation}\label{550}\sum_{j=1}^{16^{m}}\mathbb{E}[\widetilde{X}_{j}1_{\{X_{j}>\varepsilon\}}|\Phi^{S_{n},\widetilde{S}_{n,m}}]\to0\end{equation}
in probability as $m\to\infty$.
\end{lemma}
\begin{proof}
We need to deal with the maximum of the field $\Phi^{S_{n},\widetilde{S}_{n,m}}$. It follows the same line of the proof of \cite[Lemma 4.4]{ref1} that there is a constant $c=c(\delta)$ such that
\begin{equation}\label{551}\lim_{m\to\infty}\mathbb{P}\left(\sup_{x\in\widetilde{S}_{n,m}^{\delta}}\Phi^{S_{n},\widetilde{S}_{n,m}}(x)>2\sqrt{2\gamma}\log(2^{m})+c\sqrt{\log(2^{m})}\right)=0.\end{equation}
\par
Then, let $A_{n,m}$ be the event in \eqref{551} and applying the uniform bound
\begin{equation}\label{phi}\sup_{x\in\widetilde{S}_{n,m}^{\delta}}\mathrm{var}[\Phi^{S_{n},\widetilde{S}_{n,m}}(x)]\leq c+\gamma\log(2^{m})\end{equation}
on $A_{n,m}$, we have
\begin{align*}
\mathbb{E}[1_{A_{n,m}^{c}}\widetilde{X}_{j}1_{\{X_{j}>\varepsilon\}}|\Phi^{S_{n},\widetilde{S}_{n,m}}]&\leq\frac{1}{\varepsilon}\mathbb{E}[1_{A_{n,m}^{c}}\widetilde{X}_{j}^{2}|\Phi^{S_{n},\widetilde{S}_{n,m}}]\\
&\leq c\frac{\|f\|^{2}}{\varepsilon}e^{4\lambda\log(2^{m})+c(\delta)\sqrt{\log(2^{m})}}e^{8\lambda^{2}\log(2^{m})}\mathbb{E}[\widehat{Z}_{\lambda}^{S_{n+m,j,M}}(S_{n+m,j}^{\delta})^{2}].
\end{align*}
Then, Proposition \ref{4.2} and Lemma \ref{3.2} ensure that, for some $c,c'=c,c'(M,n,\delta)$,
\[\mathbb{E}[\widehat{Z}_{\lambda}^{S_{n+m,j,M}}(S_{n+m,j}^{\delta})^{2}]\leq c(\mathbb{E}[Z_{\lambda}^{S_{n+m,j}}(S_{n+m,j})])^{2}\leq c'\cdot 2^{-m(8+8\lambda^{2})}.\]
For any $\zeta>0$, we thus get
\[\mathbb{P}\left(\sum_{j=1}^{16^{m}}\mathbb{E}[\widetilde{X}_{j}1_{\{X_{j}>\varepsilon\}}|\Phi^{S_{n},\widetilde{S}_{n,m}}]\geq\zeta\right)\leq\mathbb{P}(A_{n,m})+\frac{c''}{\varepsilon\zeta}2^{-4m(1-\lambda)^{2}}e^{c(\delta)\sqrt{\log(2^{m})}}.\]
Thanks to \eqref{551}, let $m\to\infty$, we finish the proof of Lemma \ref{3.12}.
\end{proof}
\par
We now come back to the proof of ``$\geq$''. Using \eqref{542} and \eqref{550} together, we get
\begin{equation}\label{543}\limsup\limits_{m\to\infty}\mathbb{E}e^{-\langle\widetilde{Z}_{\lambda}^{S_{n},M,m},f_{m,\delta}\rangle}\leq\limsup\limits_{m\to\infty}\mathbb{E}\left[\exp\left(-e^{-\varepsilon}\mathbb{E}[\langle\widetilde{Z}_{\lambda}^{S_{n},M,m},f_{m,\delta}\rangle|\Phi^{S_{n},\widetilde{S}_{n,m}}]\right)\right].\end{equation}
We now need to prove
\begin{equation}\label{552}\mathbb{E}[\langle Z_{\lambda}^{S_{n}},f_{m,\delta}\rangle|\Phi^{S_{n},\widetilde{S}_{n,m}}]-\mathbb{E}[\langle \widetilde{Z}_{\lambda}^{S_{n},M,m},f_{m,\delta}\rangle|\Phi^{S_{n},\widetilde{S}_{n,m}}]\to0\end{equation}
in probability as $m\to\infty$ followed by $M\to\infty$. Then we will prove \eqref{552} by showing the left-hand side, after taking expectation, converges to 0.
\par
Invoking \eqref{phi}, by the property of conditional expectation, the expectation of the left-hand side of \eqref{552} can be bounded by
\[\mathbb{E}\langle Z_{\lambda}^{S_{n}}-\widehat{Z}_{\lambda}^{S_{n},M},f_{m,\delta}\rangle\leq c\cdot2^{m(4\lambda^{2}+4)}\|f\|\mathbb{E}[Z_{\lambda}^{S_{n+m}}(S_{n+m})-\widehat{Z}_{\lambda}^{S_{n+m,M}}(S_{n+m})].\]
Proposition \ref{4.1} yields
\[\mathbb{E}[Z_{\lambda}^{S_{n+m}}(S_{n+m})-\widehat{Z}_{\lambda}^{S_{n+m,M}}(S_{n+m})]\leq ce^{-\widetilde{c}M^{2}}2^{-m(4+4\lambda^{2})},\]
which proves \eqref{552}. Recall \eqref{540}, taking $M\to\infty$ and $\varepsilon\downarrow0$, using \eqref{552} in \eqref{543}, we finally get
\[\mathbb{E}e^{-\langle Z_{\lambda}^{S_{n}},f\rangle}\leq\mathbb{E}e^{-\langle Y_{m}^{S_{n}},f_{m,\delta}\rangle}.\]
Then taking $m\to\infty$ and $\delta\downarrow0$, we have
\[\mathbb{E}e^{-\langle Z_{\lambda}^{S_{n}},f\rangle}\leq\mathbb{E}e^{-\langle Y_{\infty}^{S_{n}},f\rangle}.\]
We now finish the proof of ``$\leq$''. Combining with the proof of ``$\geq$'', we complete the proof of Proposition \ref{3.11}.
\end{proof}
\par
We can finally prove Theorem \ref{2.1}:
\par
\begin{proof}[Proof of Theorem \ref{2.1}.] Proposition \ref{3.7} shows that any subsequential limit $\eta^{D}$ can be factorized as $\eta^{D}(\mathrm{d}x\mathrm{d}h)=Z_{\lambda}^{D}(\mathrm{d}x)\otimes e^{-\pi\lambda h}\mathrm{d}h$, where $Z_{\lambda}^{D}$ measures obey properties (1-6) from Proposition \ref{3.8}. And Proposition \ref{3.11} finally determined the law of $Z_{\lambda}^{D}$ for $D$ being any dyadic cube, particularly for $D=(0,1)^{4}$, which completes the proof.\
\end{proof}
\section{Characterization of the limit measure}
In this part, we will give the proof of Theorem \ref{2.4}, which shows that the limit measure $Z_{\lambda}^{D}$ has the law of the Gaussian multiplicative chaos on $D$.
\par
\begin{proof}[Proof of Theorem \ref{2.4}]
In fact, we can rewrite $Y_{m}^{S_{n}}$ defined in \eqref{ym} as
\[Y_{m}^{S_{n}}=\frac{\sqrt{\pi}}{4}e^{\frac{4\lambda^{2}}{\gamma}s_{S_{n}}(x)}\sum_{j=1}^{16^{m}}e^{\pi\lambda\Phi^{S_{n},\widetilde{S}_{n,m}(x)}-\frac{1}{2}\pi^{2}\lambda^{2}\mathbb{E}[\Phi^{S_{n},\widetilde{S}_{n,m}}(x)^{2}]}1_{S_{n+m,j}(x)}\mathrm{d}x.\]
Recall that $\Phi^{S_{n},\widetilde{S}_{n,m}}$ is the sum of independent fields:
\[\Phi^{S_{n},\widetilde{S}_{n,m}}=\sum_{j=1}^{m}\Phi^{\widetilde{S}_{n,j-1},\widetilde{S}_{n,j}}.\]
From now on, for brevity of notation, we will write $S$ (resp. $S^{j}$) instead of $S_{n}$ (resp. $\widetilde{S}_{n,j}$). Now, let $H_{k}$ denote the subspace of $H_{0}^{2}(S)$, which contains the functions that are biharmonic in $S^{k}$, with Dirichlet boundary conditions on $\partial S^{k-1}$. We claim that:
\[H_{0}^{2}(S)=\bigoplus_{k=0}^{\infty}H_{k}.\]
First of all, we need to show that $H_{i}$ and $H_{j}$ are orthogonal to each other for all $i\neq j$. Invoking Gauss-Green formula, for $f\in H_{i},g\in H_{j}(i<j)$, we have:
\[\langle f,g\rangle_{\Delta}=\int_{S}\Delta f(x)\Delta g(x)\mathrm{d}x=\int_{S^{i}}g\Delta^{2}f\mathrm{d}x+\int_{\partial S^{i}}\left(\Delta f\frac{\partial g}{\partial\bm{n}}-g\frac{\partial(\Delta f)}{\partial\bm{n}}\right)\mathrm{d}\sigma=0\]
since $\Delta^{2}f=0$ in $S^{i}$ and $g=\frac{\partial g}{\partial\bm{n}}=0$ on $\partial S^{i}$, thus we prove the orthogonality. Then we need to show that $H_{0}^{2}(S)\backslash\bigoplus_{k=0}^{\infty}H_{k}$ is trivial. Assuming $f$ is in the orthocomplement space of $\bigoplus_{k=0}^{\infty}H_{k}$, then Gauss-Green formula implies that $f$ is the weak solution of the following boundary value problem:
\begin{equation}\label{611}\begin{cases}\Delta^{2}f(x)=0,\ \ & x\in S^{\infty},\\ f(x)=0,\ \ & x\in\partial S^{\infty},\\ \frac{\partial f}{\partial\bm{n}}(x)=0,\ \ & x\in\partial S^{\infty},\end{cases}\end{equation}
where $S^{\infty}:=\bigcap\limits_{k=0}^{\infty}S^{k}$. Thanks to the uniqueness of the solution of \eqref{611}, we conclude that $f=0$ almost everywhere, and thus complete the proof of the claim above.
\par
Then, let $\{\varphi_{k,j}\}_{j}$ denote the orthonormal basis of $H_{k-1}$, a straightforward calculation of covariance shows that:
\[\Phi^{S^{k-1},S^{k}}\overset{\mathrm{law}}{=}\sum_{j}\varphi_{k,j}Z_{k,j}\]
with $Z_{k,j}$ i.i.d. standard normals. Then we get:
\[\Phi^{S,S^{i}}\overset{\mathrm{law}}{=}\sum_{k=0}^{i-1}\sum_{j}\varphi_{k,j}Z_{k,j}.\]
Denoting $\mathscr{F}_{k}:=\sigma(Z_{i,j}:i+j\leq k)$, we immediately get:
\[\mathbb{E}[Y_{m}^{S}(A)|\mathscr{F}_{k}]=\int_{A}\frac{\sqrt{\pi}}{4}e^{\frac{4\lambda^{2}}{\gamma}s_{S}(x)}\mu_{k(k-1)/2}^{S,\pi\lambda}(\mathrm{d}x),\quad m\geq k,\]
where $\mu_{n}^{S,\beta}$ is the Gaussian multiplicative chaos measure defined through the orthonormal basis $\{\varphi_{k,j}\}_{k,j}$ rearranged through the complete order:
\[(k,j)\preceq(k',j')\Leftrightarrow k<k' \;\mathrm{or}\; k=k',j<j'.\]
Letting $m\to\infty$, since we have proved $Y_{m}^{S}(A)\to Y_{\infty}^{S}(A)$ a.s. in Lemma \ref{3.10}, we get:
\[\mathbb{E}[Y_{\infty}^{S}(A)|\mathscr{F}_{k}]=\int_{A}\frac{\sqrt{\pi}}{4}e^{\frac{4\lambda^{2}}{\gamma}s_{S}(x)}\mu_{k(k-1)/2}^{S,\pi\lambda}(\mathrm{d}x).\]
Then, letting $k\to\infty$ on both sides, by the triviality of $\mathscr{F}_{\infty}:=\sigma(\bigcup_{k\geq1}\mathscr{F}_{k})$, we get:
\[Y_{\infty}^{S}(A)=\int_{A}\frac{\sqrt{\pi}}{4}e^{\frac{4\lambda^{2}}{\gamma}s_{S}(x)}\mu_{\infty}^{S,\pi\lambda}(\mathrm{d}x).\]
By the arbitrariness of $A$, we complete the proof since $Y_{\infty}^{S}$ and $Z_{\lambda}^{S}$ has the same law.
\end{proof}
\appendix
\section{The Gibbs-Markov property}
In this Appendix, we will give a detailed proof of the Gibbs-Markov property of the membrane model (Theorem \ref{gmp}). For a given finite set $V\subset\mathbb{Z}^{d}$, consider the Hilbert space $\mathcal{H}^{V}:=\{f:\mathbb{Z}^{d}\to\mathbb{R}:\mathrm{supp}(f)\subset V\}$ endowed with the inner product
\[\langle f,g\rangle_{\Delta}:=\sum_{x\in\mathbb{Z}^{d}}\Delta f(x)\cdot\Delta g(x).\]
We then have the following lemma:
\begin{lemma}\label{A.1}
For the setting as above with $V$ finite, let $\{\varphi_{n}:n=1,2,\cdots,|V|\}$ be an orthonormal basis in $\mathcal{H}^{V}$ and let $Z_{1},\cdots,Z_{|V|}$ be i.i.d. standard normals. Then
\begin{equation}\label{a11}\widetilde{h}_{x}^{V}:=\sum_{n=1}^{|V|}\varphi_{n}(x)Z_{n},\quad x\in\mathbb{Z}^{d}\end{equation}
has the law of the membrane model in $V$.
\end{lemma}
\begin{proof}
Since $\{\widetilde{h}_{x}^{V}\}_{x\in\mathbb{Z}^{d}}$ is a centered Gaussian vanishing outside $V$ and
\[\mathbb{E}[\widetilde{h}_{x}^{V}\widetilde{h}_{y}^{V}]=\sum_{n=1}^{|V|}\varphi_{n}(x)\varphi_{n}(y),\]
it suffices to prove that
\[G^{V}(x,y)=\sum_{n=1}^{|V|}\varphi_{n}(x)\varphi_{n}(y).\]
Fixing $y$ on the left-hand side, the right-hand side can be seen as the Fourier expansion of $G^{V}(x,y)$. Assume that
\[G^{V}(x,y)=\sum_{n=1}^{|V|}c_{n}\varphi_{n}(x),\]
then
\[c_{n}=\sum_{x\in\mathbb{Z}^{d}}\Delta G^{V}(x,y)\cdot\Delta\varphi_{n}(x)=\sum_{x\in\mathbb{Z}^{d}}\Delta^{2}G^{V}(x,y)\cdot\varphi_{n}(x)=\sum_{x\in\mathbb{Z}^{d}}\delta_{y}(x)\varphi_{n}(x)=\varphi_{n}(y).\]
Since $y$ is arbitrary, we complete the proof of the lemma.
\end{proof}
\par
Then we will prove a decomposition lemma, which leads to the proof of Theorem \ref{gmp}.
\begin{lemma}\label{A.2}
Let $U\subset V\subset\subset\mathbb{Z}^{d}$, and $\mathcal{H}^{V},\mathcal{H}^{U}$ as above. Let $\mathrm{Bih}(U)$ be the space of discrete biharmonic functions in $U$, then we have:
\[\mathcal{H}^{V}=\mathcal{H}^{U}\oplus\mathrm{Bih}(U).\]
\end{lemma}
\begin{proof}
First observe that any functions $f\in\mathcal{H}^{U}$ is orthogonal to any functions in $\mathrm{Bih}(U)$ by the Gauss-Green formula in the discrete case. Then for any $f\in\mathcal{H}^{V}$, letting $f_{0}$ be the orthogonal projection of $f$ onto $\mathcal{H}^{U}$, it suffices to show that $\varphi=f-f_{0}$ is discrete biharmonic on $U$. Note that $\varphi$ is orthogonal to $\mathcal{H}^{U}$ by definition, then for any test function $\psi\in\mathcal{H}^{U}$, by the Gauss-Green formula in the discrete case, we have
\[\sum_{x\in U}(\Delta^{2}\varphi(x))\psi(x)=\sum_{x\in\mathbb{Z}^{d}}\Delta\varphi(x)\cdot\Delta\psi(x)=0.\]
Since $\psi$ is arbitrary, $\varphi$ is discrete biharmonic in $U$. We thus complete the proof.
\end{proof}
\par
By the two lemmas above, we can now prove Theorem \ref{gmp}.
\par
\begin{proof}[Proof of Theorem \ref{gmp}.]
Let $\{f_{n}\}$ be an orthonormal basis of $\mathcal{H}^{U}$, and $\{\phi_{n}\}$ be an orthonormal basis of $\mathrm{Bih}(U)$. We can write
\[\widetilde{h}^{U}:=\sum X_{n}f_{n},\quad\varphi^{V,U}:=\sum Y_{n}\phi_{n},\]
where $X_{n},Y_{n}$ are i.i.d. standard Gaussian variables. Then by Lemma \ref{A.1}, $\widetilde{h}^{U}$ has the law of the membrane model in $U$. Thanks to Lemma \ref{A.2}, $\{f_{n},\phi_{n}\}$ is an orthonormal basis of $\mathcal{H}^{V}$, then again invoking Lemma \ref{A.1}, $\widetilde{h}^{U}+\varphi^{V,U}$ has the law of the membrane model in $V$, since $X_{n}$ and $Y_{n}$ are independent, we get the independence between $\widetilde{h}^{U}$ and $\varphi^{V,U}$. Since $\phi_{n}\in\mathrm{Bih}(U)$, we prove that $\varphi^{V,U}$ is discrete biharmonic in $U$.
\par
Note that $\varphi^{V,U}$ has the law of the biharmonic extension of $h^{V}(x)(x\in V\backslash U)$ in $U$ since $\phi_{n}\in\mathrm{Bih}(U)$ and $h^{U}$ vanishes outside $U$. In order to complete the proof, we only need to show that $\mathbb{E}[h^{V}|\sigma(h_{z}^{V}:z\in V\backslash U)]$ is discrete biharmonic on $U$. Fix $x\in U$, recall Definition \ref{dmm}, since
\[\sum_{y\in\mathbb{Z}^{d}}|\Delta h_{y}|^{2}=(1+1/2d)h_{x}^{2}+h_{x}\left(-\frac{1}{d}\sum_{i}h_{x\pm e_{i}}+\frac{1}{d^{2}}\sum_{i\neq j}h_{x\pm e_{i}\pm e_{j}}+\frac{1}{2d^{2}}\sum_{i}h_{x\pm2e_{i}}\right)+\cdots,\]
where $h_{x}$ does not occur in $\cdots$ terms. A simple calculation shows that:
\[\mathbb{E}[h_{x}^{V}|\sigma(h_{z}^{V}:z\in V\backslash\{x\})]=\frac{d}{2d+1}\left(\frac{1}{d}\sum_{i}h_{x\pm e_{i}}-\frac{1}{d^{2}}\sum_{i<j}h_{x\pm e_{i}\pm e_{j}}-\frac{1}{2d^{2}}\sum_{i}h_{x\pm2e_{i}}\right).\]
Thanks to the ``smaller $\sigma$-algebra always wins'' property of conditional expectation, we have
\begin{align*}
\frac{2d+1}{2d}\mathbb{E}[h^{V}|\sigma(h_{z}^{V}:z\in V\backslash U)]&=\mathbb{E}[\mathbb{E}[h_{x}^{V}|\sigma(h_{z}^{V}:z\in V\backslash\{x\})]|\sigma(h_{z}^{V}:z\in V\backslash U)]\\
&=\mathbb{E}\Big[\frac{1}{2d}\sum_{i}h_{x\pm e_{i}}-\frac{1}{2d^{2}}\sum_{i<j}h_{x\pm e_{i}\pm e_{j}}-\frac{1}{4d^{2}}\sum_{i}h_{x\pm2e_{i}}\Big|\sigma(h_{z}^{V}:z\in V\backslash U)\Big].
\end{align*}
By the explicit form of the discrete bi-Laplacian operator, and the arbitrariness of $x\in U$, we prove that $\mathbb{E}[h^{V}|\sigma(h_{z}^{V}:z\in V\backslash U)]$ is discrete biharmonic on $U$, and thus complete the proof.
\end{proof}

\noindent {\bf Acknowledgments}\quad The authors thank Hao Ge and Jiaxi Zhao for helpful suggestions during various stages of this work. The authors also thank an anonymous referee for various helpful and detailed comments.

\end{document}